\documentclass{gtart}
\usepackage{amsmath,amssymb}

\usepackage{graphicx}

\theoremstyle{plain}
\newtheorem{lem}{Lemma}
\newtheorem{thm}{Theorem}
\newtheorem{cor}{Corollary}

\newcommand{\degree}{^\circ}
\newcommand{\disjoint}{\amalg}  
\newcommand{\equivalent}{\sim}
\newcommand{\other}{\overline}
\newcommand{\ttimes}{\twisted\times}
\newcommand{\twisted}{\widetilde}

\newcommand{\rZ}{\mathbb Z}

\DeclareMathOperator{\writhe}{w}
\newcommand{\wm}{\varepsilon}   
\newcommand{\p}{\hspace*{-3pt}+}
\newcommand{\n}{\hspace*{-3pt}-}

\newcommand{\boundary}{\partial}
\newcommand{\inv}[2][1]{{#2}^{-#1}} 

\begin{document}

\title{Twisted Link Theory}
\author{Mario O. Bourgoin}
\address{Department of Mathematics \\ Brandeis University \\ Waltham, MA 02454}
\email{mob@brandeis.edu}
\date{\today}
\begin{abstract}
  We introduce stable equivalence classes of oriented links in
  orientable three-manifolds that are orientation $I$-bundles over
  closed but not necessarily orientable surfaces.  We call these
  twisted links, and show that they subsume the virtual knots
  introduced by L.~Kauffman, and the projective links introduced by
  Yu.~Drobotukhina.  We show that these links have unique minimal
  genus three-manifolds.  We use link diagrams to define an extension
  of the Jones polynomial for these links, and show that this
  polynomial fails to distinguish two-colorable links over
  non-orientable surfaces from non-two-colorable virtual links.
\end{abstract}
\makeshorttitle

\section{Introduction}

\subsection{Virtual Links}

\begin{figure}[htbp]
  \centering
  \raisebox{-.75in}{\includegraphics[height=1.5in]{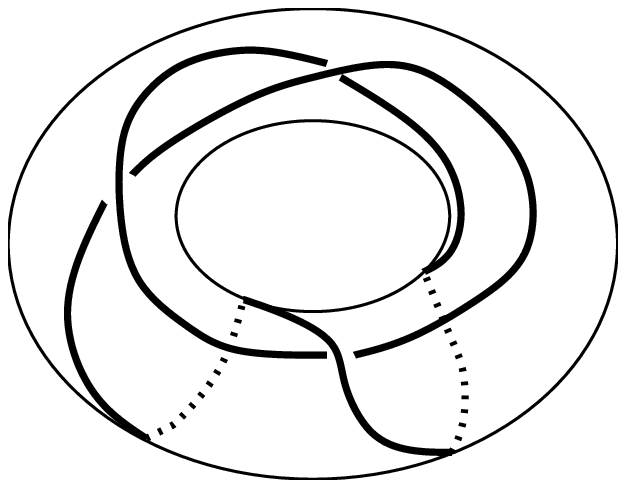}}
  $\Rightarrow$
  \raisebox{-.75in}{\includegraphics[height=1.5in]{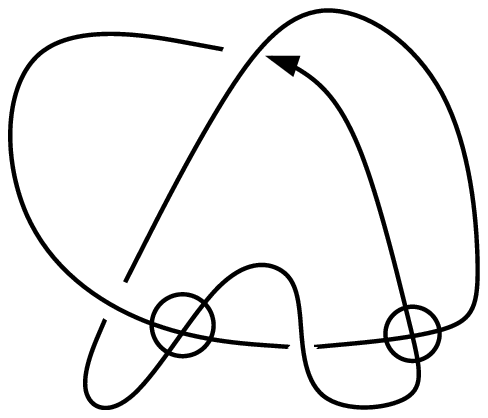}}
  \caption{A knot in a thickened torus and its diagram.}
  \label{fig:vrefoil}
\end{figure}

In 1996, Louis Kauffman introduced a generalization of classical links
to stable embeddings of a disjoint union of circles in a thickened
compact oriented surface by using the notion of an oriented Gauss
code~\cite{MR2000i:57011}.  An example of such a link is shown in
Figure~\ref{fig:vrefoil}.
Kauffman was motivated in part by the desire to allow all oriented
Gauss codes to be associated with diagrams of links.  When the surface
being thickened is a sphere, the links are classical, however there
are virtual links that are not classical links.  Since by a result of
Kuperberg~\cite{MR1997331}, virtual links have unique irreducible
representatives, virtual link theory is a proper extension of
classical links.

Many invariants of classical links formally extend to virtual links
through the use of their diagrams by ignoring the virtual crossings,
although the resulting invariant may have different properties than
those of the classical version.  A virtual link can be formally
associated with a group through a Wirtinger presentation obtained from
any diagram for the link.  However, the group of a virtual link may
not be residually finite, so it may not be the fundamental group of
any three-manifold~\cite{MR2002m:57011}.  A virtual link has a Jones
polynomial, but the classical relation between the polynomial's
exponents and the number of link components only holds if the link has
a two-colorable diagram as was shown by N.~Kamada~\cite{MR1914297}.

\subsection{Projective Links}

\begin{figure}[htbp]
  \centering
  \raisebox{-.75in}{\includegraphics[height=1.5in]{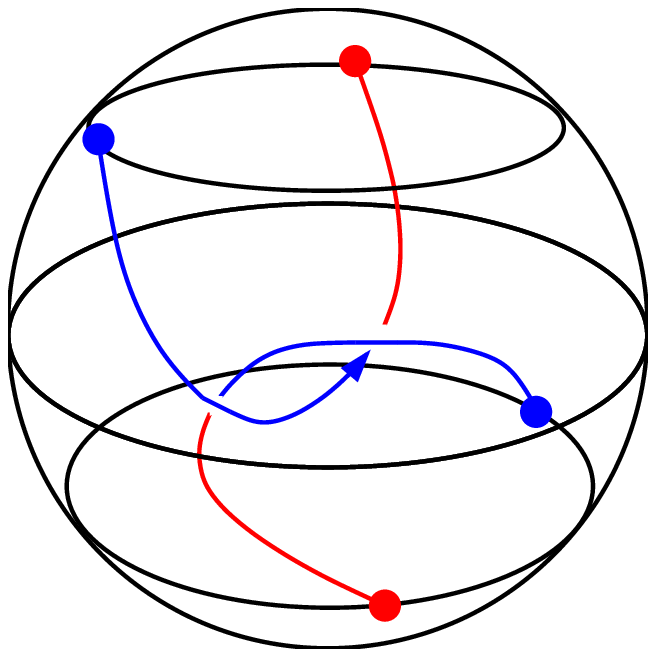}}
  $\Rightarrow$
  \raisebox{-.75in}{\includegraphics[height=1.5in]{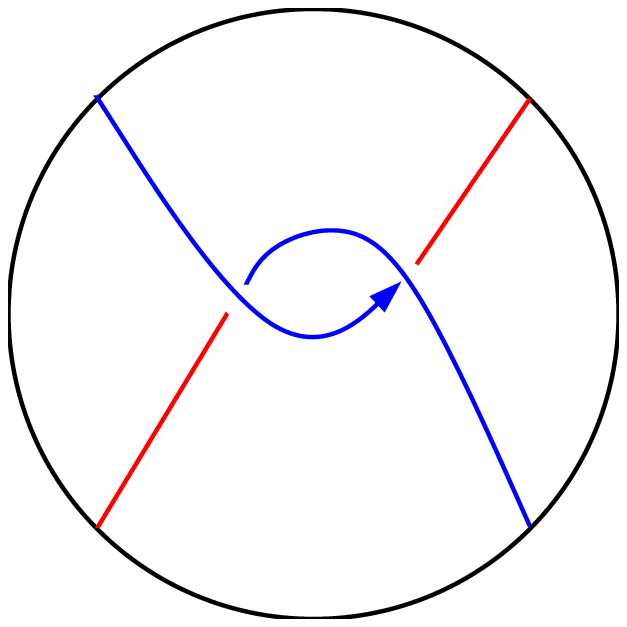}}
  \caption{A projective knot and its diagram.}
  \label{fig:twofoil}
\end{figure}

In 1990, Yu.~V.~Drobotukhina introduced the study of links in real
projective space as a generalization of links in the
three-sphere~\cite{MR91i:57001}.  On the left of
Figure~\ref{fig:twofoil},
a non-trivial projective link is shown using the 3-ball model of
projective space.  Drobotukhina showed how to create diagrams of
projective links in 2-disks, and she extended Reidemeister moves on
planar diagrams to include moves across the boundary of the disks.  On
the right of Figure~\ref{fig:twofoil}, is shown a diagram of the link
in the projective plane.  Drobotukhina showed that projective links
admit a Jones polynomial invariant and studied some of its properties.

\subsection{Links in Oriented Thickenings}

\begin{figure}[htbp]
  \centering
  \raisebox{-.75in}{\includegraphics[height=1.5in]{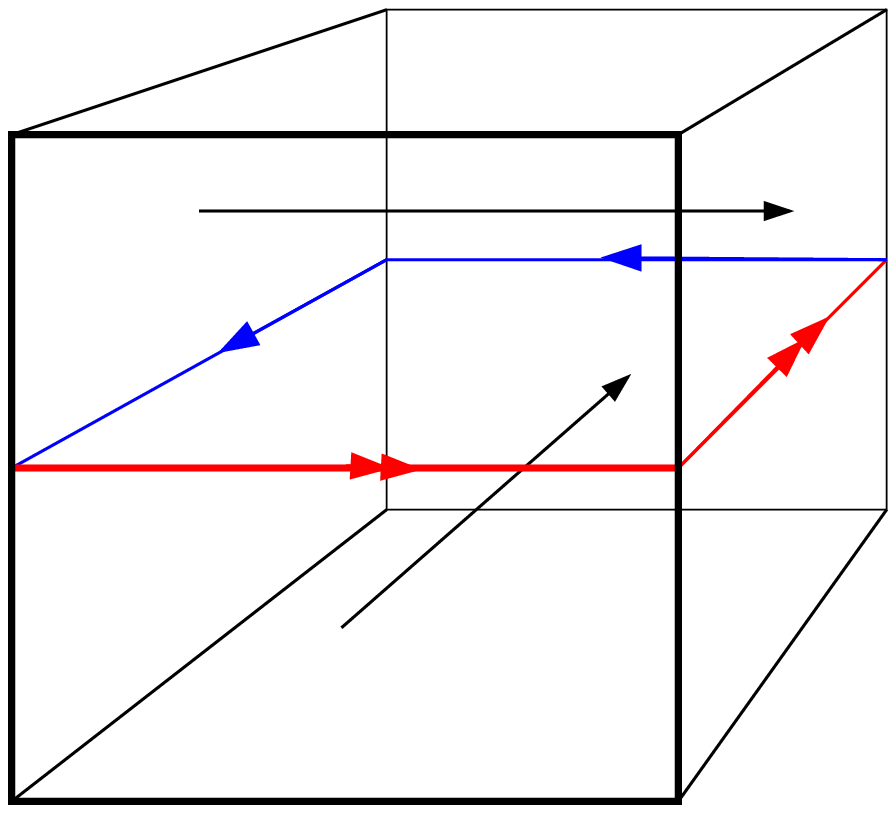}}
  $\Longrightarrow$
  \raisebox{-.75in}{\includegraphics[height=1.5in]{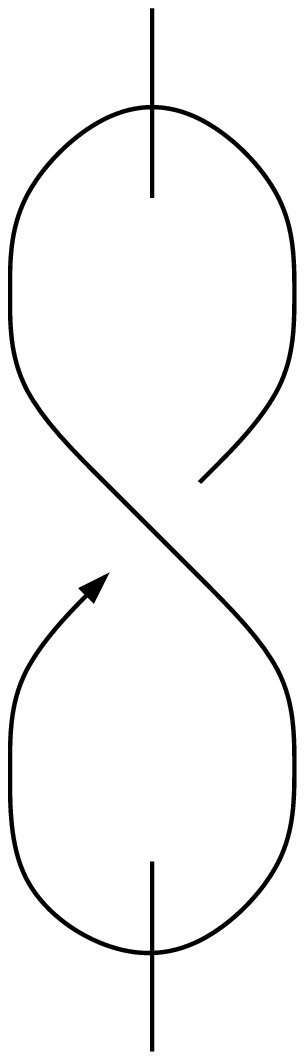}}
  \caption{A onefoil knot in a thickened Klein bottle.}
  \label{fig:onefoil}
\end{figure}

We introduce \emph{links in oriented thickenings} as stable ambient
isotopy classes of oriented circles in oriented three-manifolds that
are orientation $I$-bundles over closed but not necessarily orientable
or connected surfaces.  Figure~\ref{fig:onefoil}
shows a onefoil knot in a thickened Klein bottle on the left.  The
thickened Klein bottle is shown as a cube with identification of its
sides with the double arrows with a $180\degree$ turn, and likewise
for the sides with the single arrow.  A diagram for the onefoil is
shown on the right of the figure, with the bars on the edges showing
that paths around those edges are orientation-reversing in a
projection of the link on the embedded Klein bottle.

By considering destabilization of the oriented thickening along annuli
and M\"obius bands in the complement of the link, we get an extension
to links in oriented thickenings of a result of Greg
Kuperberg~\cite{MR1997331} for virtual links.

\begin{thm}\label{thm:unique_destabilization}
  Links in oriented thickenings have a unique irreducible
  representative.
\end{thm}

Thus we may speak of the minimum Euler genus\footnote{The Euler genus
  of a surface is defined as two minus its Euler characteristic.} of a
link in oriented thickening.  The following is an immediate corollary:

\begin{cor}
  Links in their minimum Euler genus oriented thickenings that are
  equivalent through stabilizations are equivalent without
  stabilizations.
\end{cor}

In particular, classical, projective, and virtual link theories inject
in the theory of link in oriented thickenings.

\begin{figure}[htbp]
  \raggedright\textbf{Classical Moves:}\\[1ex] \centering
  \raisebox{-.35in}{\includegraphics[height=.8in]{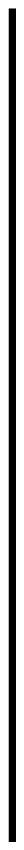}}    
  $\buildrel\text{R1}\over\Longleftrightarrow$
  \raisebox{-.35in}{\includegraphics[height=.8in]{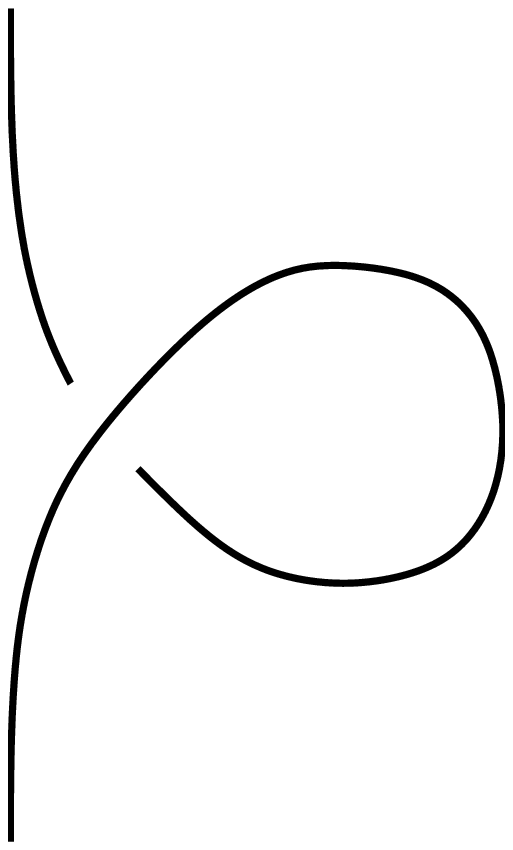}}
  \hspace{1ex}
  \raisebox{-.35in}{\includegraphics[height=.8in]{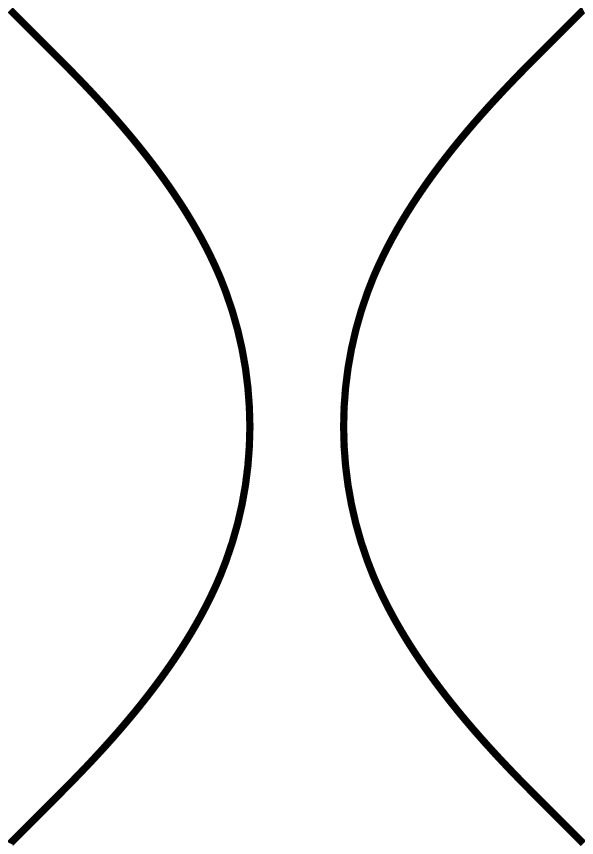}}
  $\buildrel\text{R2}\over\Longleftrightarrow$
  \raisebox{-.35in}{\includegraphics[height=.8in]{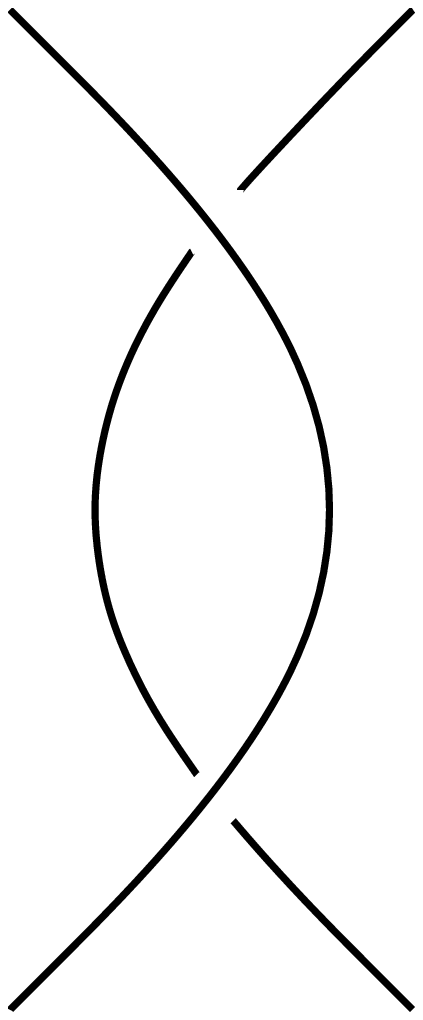}}
  \hspace{1ex}
  \raisebox{-.35in}{\includegraphics[height=.8in]{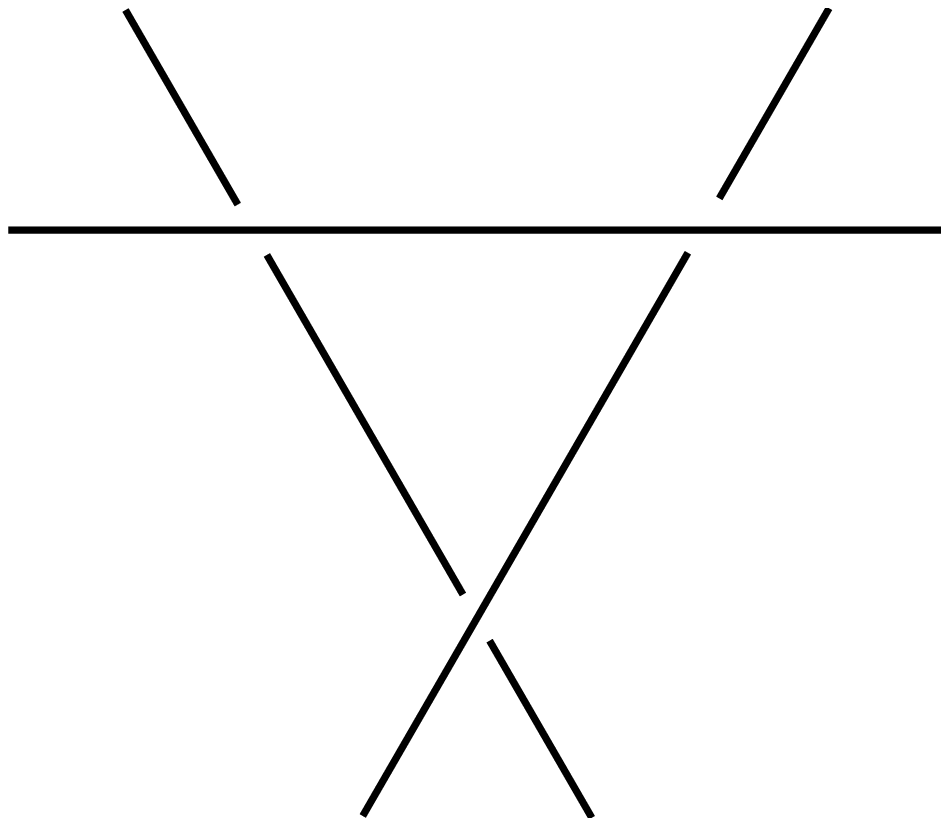}}
  $\buildrel\text{R3}\over\Longleftrightarrow$
  \raisebox{-.35in}{\includegraphics[height=.8in,angle=180,origin=c]{3lr}} \\[1ex]
  \raggedright\textbf{Virtual Extension:}\\[1ex] \centering
  \raisebox{-.35in}{\includegraphics[height=.8in]{1vert}}   
  $\buildrel\text{V1}\over\Longleftrightarrow$
  \raisebox{-.35in}{\includegraphics[height=.8in]{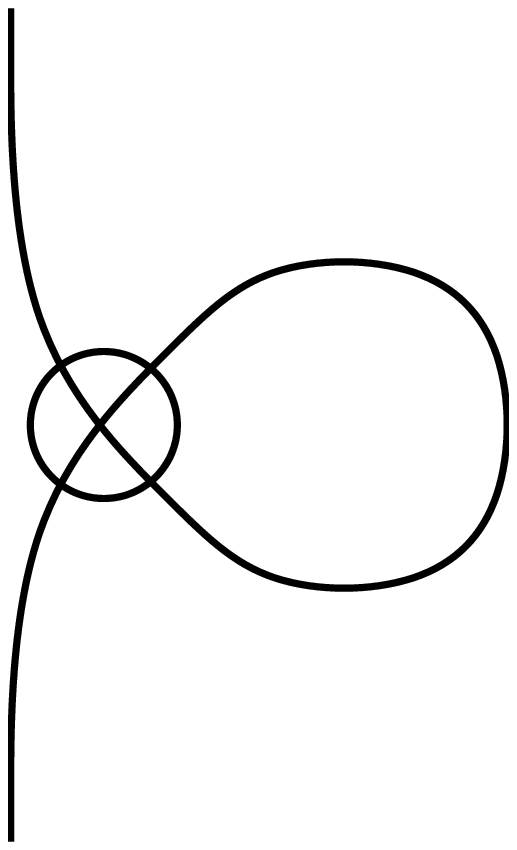}}
  \hspace{1ex}
  \raisebox{-.35in}{\includegraphics[height=.8in]{2vert}}
  $\buildrel\text{V2}\over\Longleftrightarrow$
  \raisebox{-.35in}{\includegraphics[height=.8in]{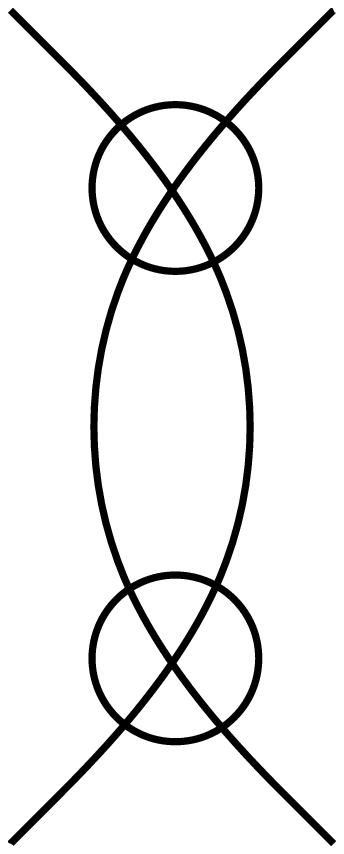}}
  \hspace{1ex}
  \raisebox{-.35in}{\includegraphics[height=.8in]{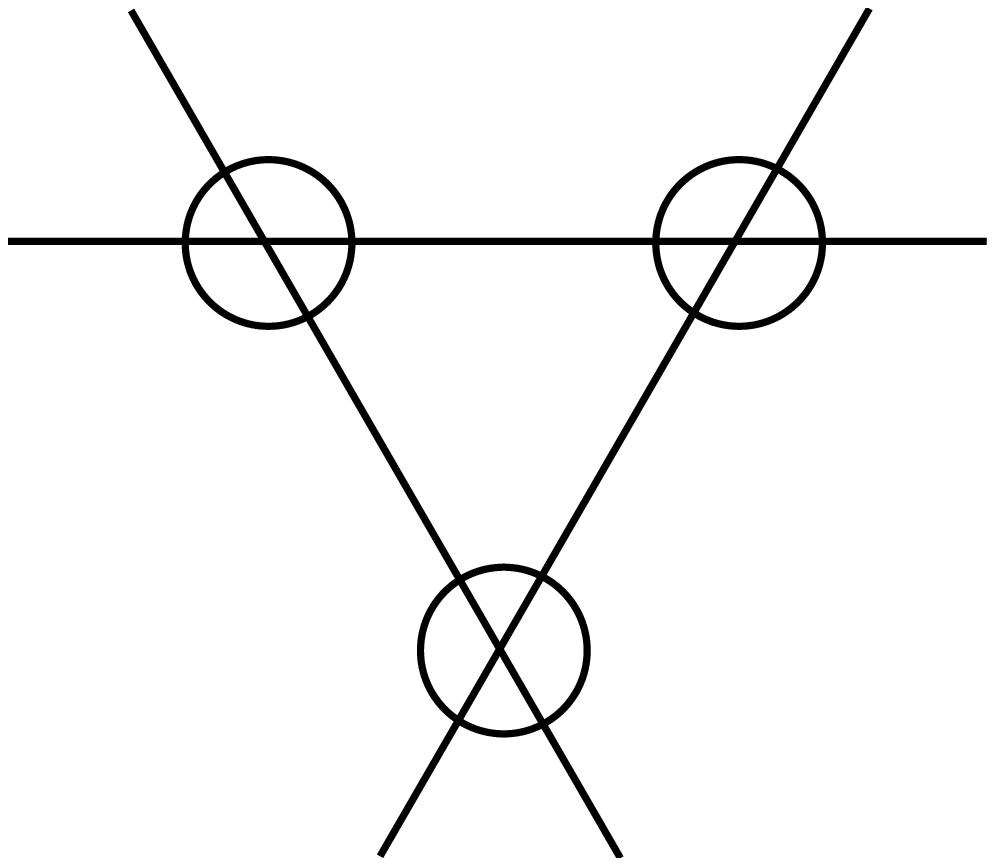}}
  $\buildrel\text{V3}\over\Longleftrightarrow$
  \raisebox{-.35in}{\includegraphics[height=.8in,angle=180,origin=c]{3virt}}
  \hspace{1ex}
  \raisebox{-.35in}{\includegraphics[height=.8in]{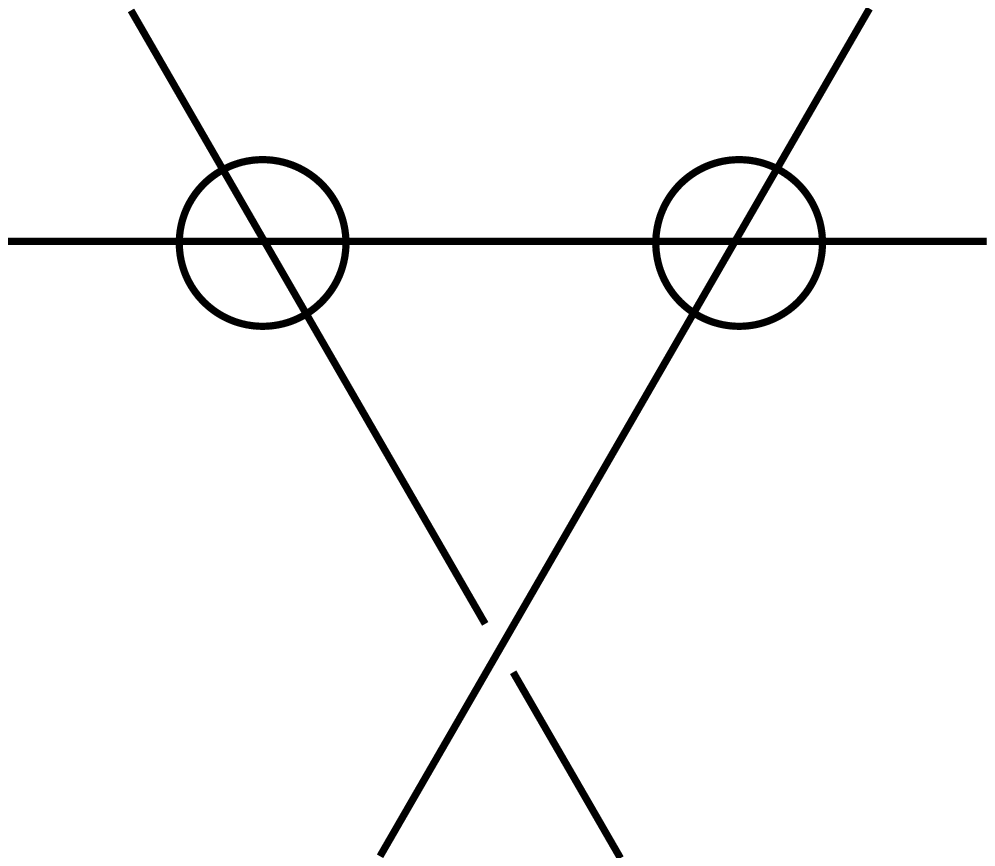}}
  $\buildrel\text{V4}\over\Longleftrightarrow$
  \raisebox{-.35in}{\includegraphics[height=.8in,angle=180,origin=c]{1lr2virt}} \\[1ex]
  \raggedright\textbf{Twisted Extension:}\\[1ex] \centering
  \raisebox{-.35in}{\includegraphics[height=.8in]{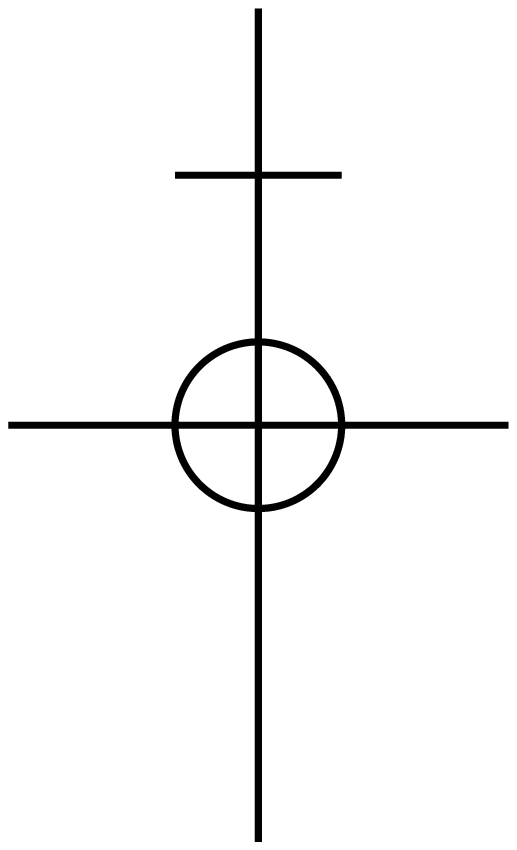}}   
  $\buildrel\text{T1}\over\Longleftrightarrow$
  \raisebox{-.35in}{\includegraphics[height=.8in,angle=180,origin=c]{1vert1virt1bar}}
  \hspace{1em}
  \raisebox{-.35in}{\includegraphics[height=.8in]{1vert}}
  $\buildrel\text{T2}\over\Longleftrightarrow$
  \raisebox{-.35in}{\includegraphics[height=.8in]{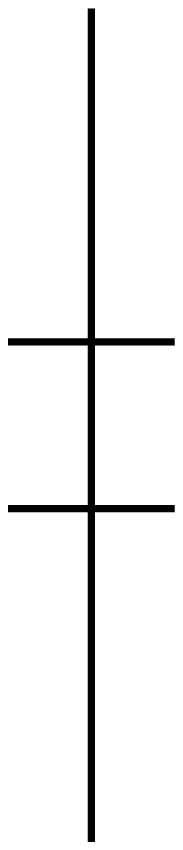}}
  \hspace{1em}
  \raisebox{-.35in}{\includegraphics[height=.8in]{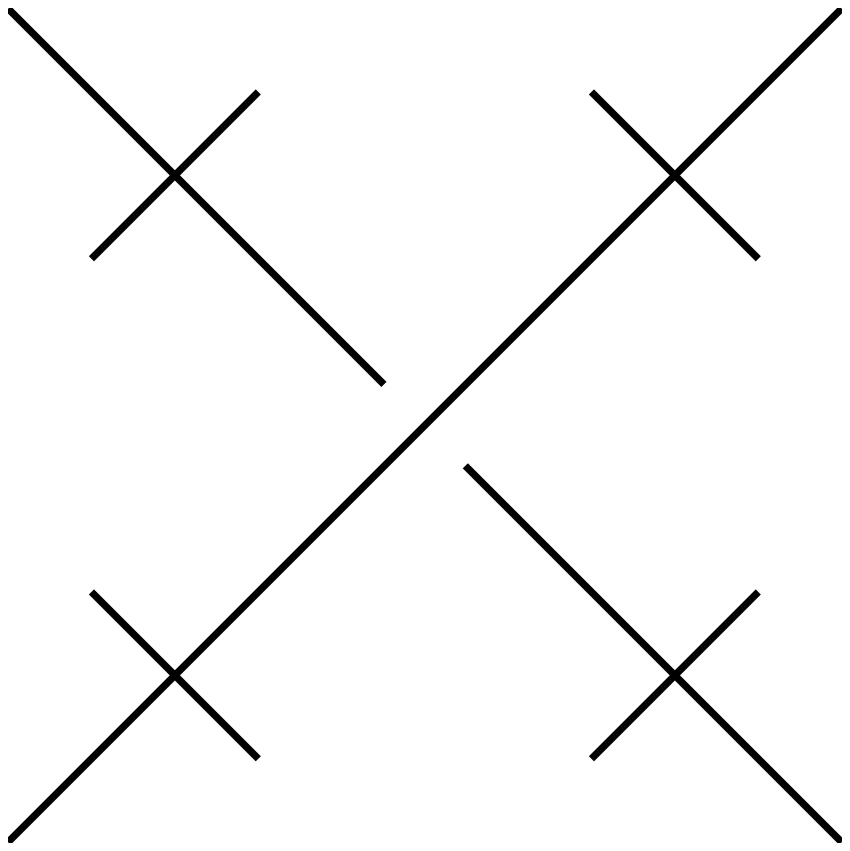}}
  $\buildrel\text{T3}\over\Longleftrightarrow$
  \raisebox{-.35in}{\includegraphics[height=.8in]{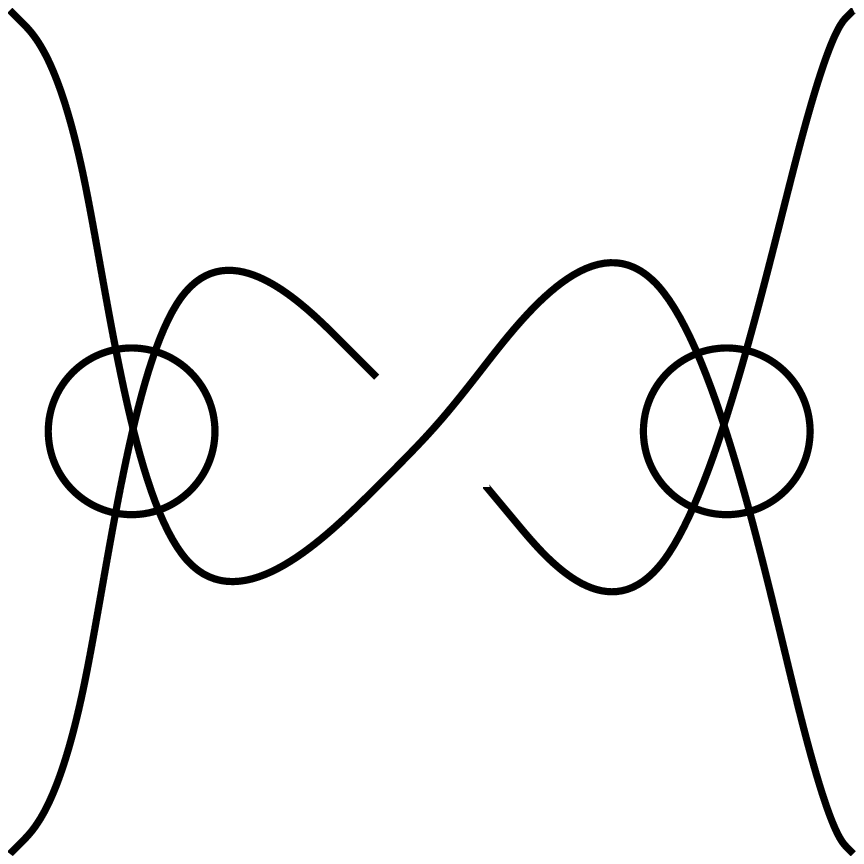}}
  \caption{The ten extended Reidemeister moves.}
  \label{fig:Reidemeister}
\end{figure}

\begin{figure}[htbp]
  \centering
  \raisebox{-.35in}{\includegraphics[height=.8in]{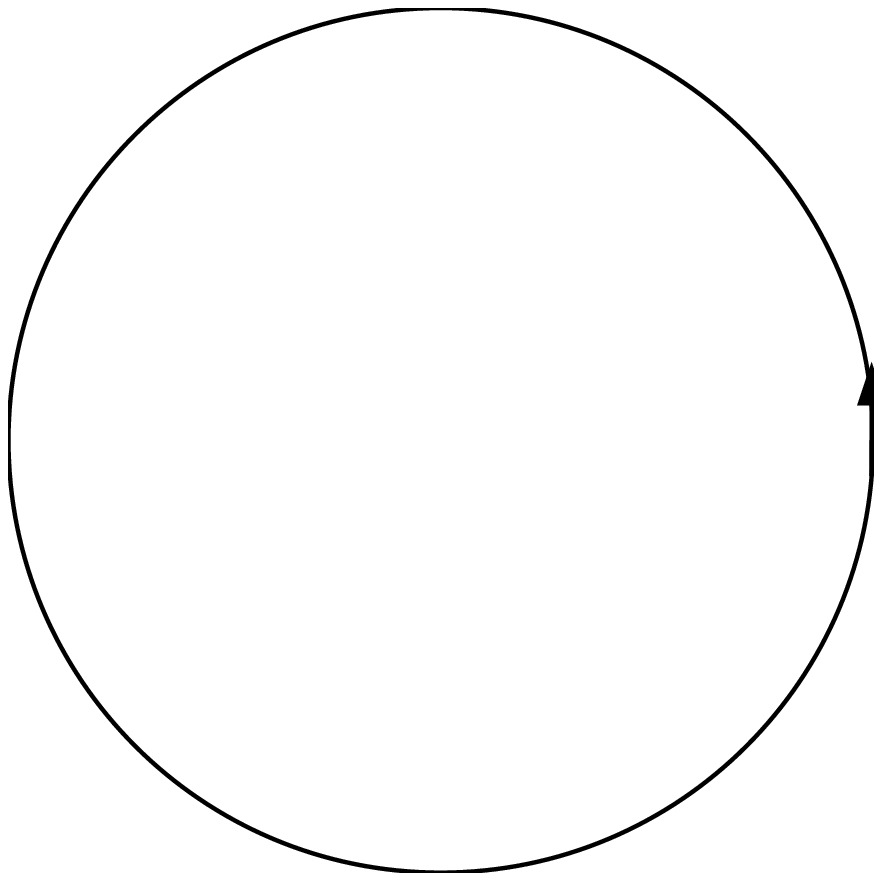}}
  $\buildrel\text{T2}\over\Longleftrightarrow$
  \raisebox{-.45in}{\includegraphics[height=1in]{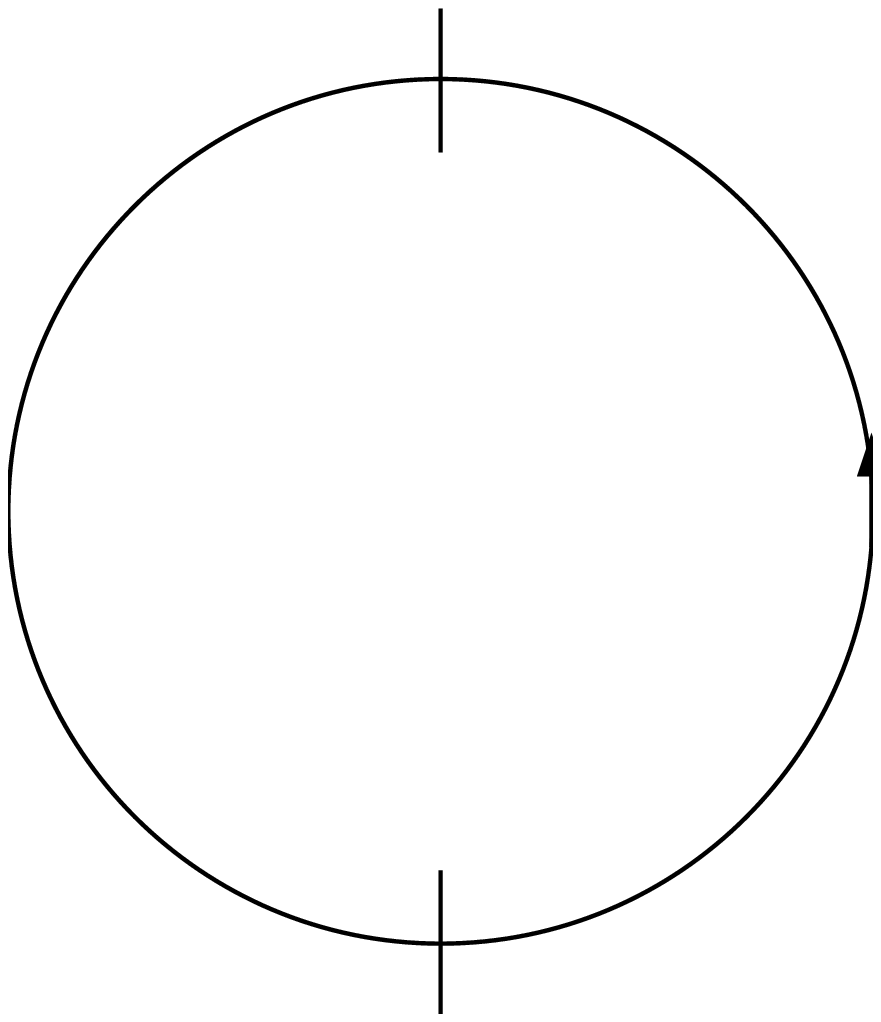}}
  $\buildrel\text{V1}\over\Longleftrightarrow$
  \raisebox{-.45in}{\includegraphics[height=1in]{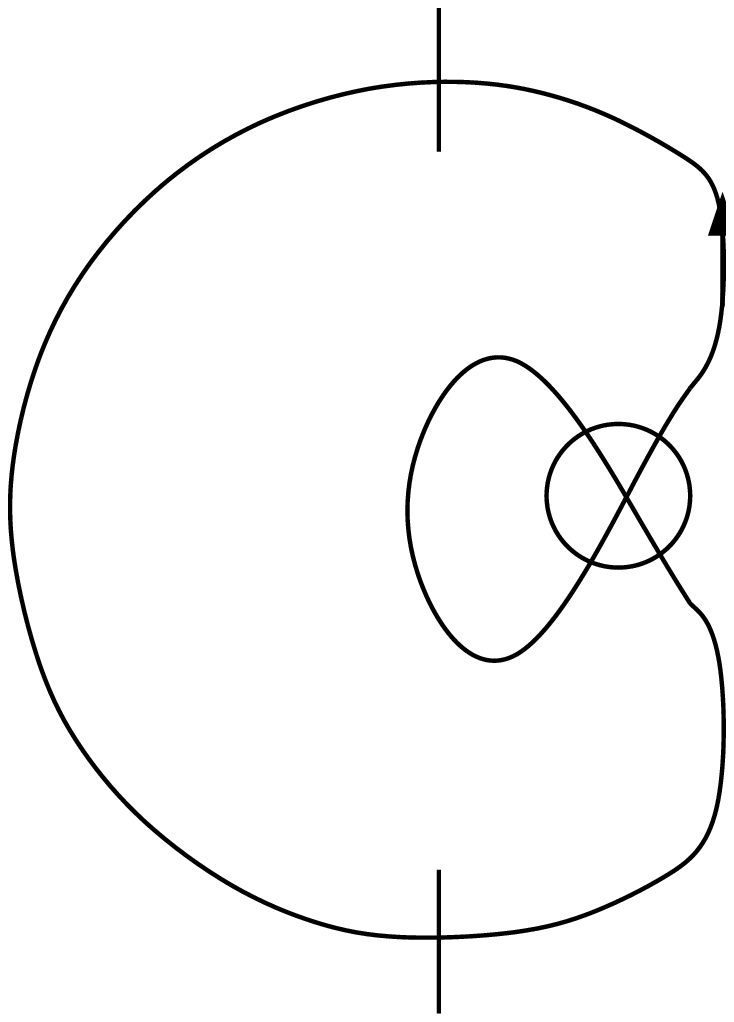}}
  $\buildrel\text{R2}\over\Longleftrightarrow$
  \raisebox{-.45in}{\includegraphics[height=1in]{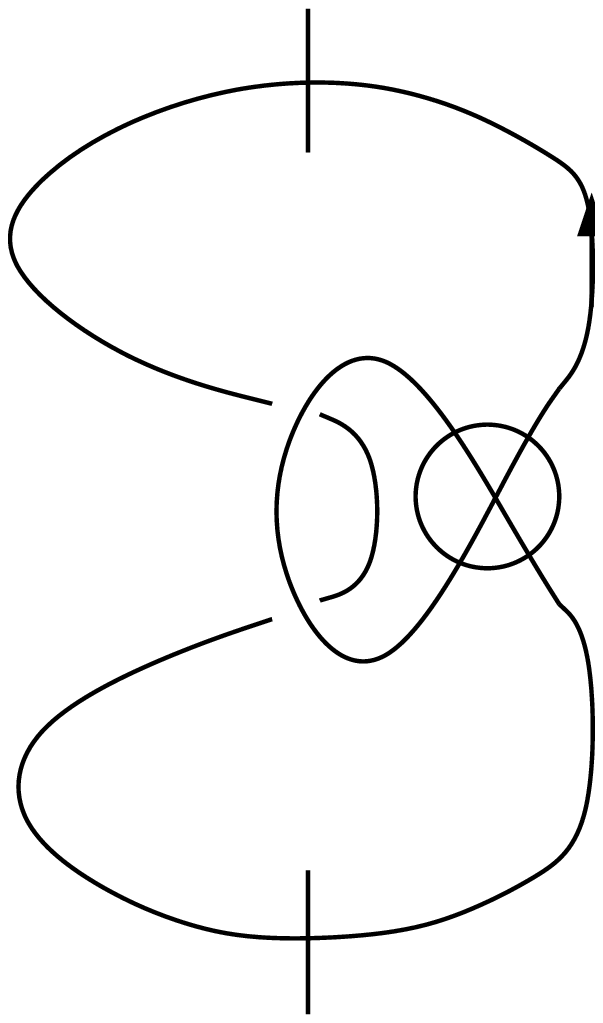}}
  \caption{Three extended Reidemeister moves on an unknot.}
  \label{fig:Ex1}
\end{figure}

We define \emph{twisted link diagrams} as marked generic planar
curves, where the markings identify the usual classical crossings,
Kauffman's virtual crossings, and now bars on edges.  We extend the
Kauffman-Reidemeister moves for virtual links to account for the new
markings by including the twisted moves described in
Figure~\ref{fig:Reidemeister}.
Since there are three classical Reidemeister moves labeled R1 to R3
and four virtual Reidemeister moves labeled V1 to V4, the three
twisted Reidemeister moves labeled T1 to T3 brings the total number of
moves to ten.  Figure~\ref{fig:Ex1}
shows an example of transforming the classical diagram of an unknot
through each of the three classes of the extended Reidemeister moves.
We will see that the leftmost diagram is for an unknot in a thickened
sphere while the rightmost diagram is for a curve in a thickened Klein
bottle.

We have the following:

\begin{thm}\label{thm:link_diagrams}
  Links in oriented thickenings correspond to classes of twisted link
  diagrams given by ambient isotopy of the diagram along with extended
  Kauffman-Reidemeister moves.
\end{thm}

Therefore, we will refer to links in oriented thickenings as
\emph{twisted links}.

\subsection{The Twisted Jones Polynomial}
\label{sec:twisted_jones}

We use a state sum to define the \emph{twisted Jones polynomial} of a
twisted link as an element of $\rZ[A^{\pm1},M]$, where the $M$
variable counts the number of circles with an odd number of bars in a
given state of the diagram.  We have the following:

\begin{thm}\label{thm:Jones_polynomial}
  The twisted Jones polynomial is an invariant of twisted links.  If a
  twisted link diagram is that of a virtual link, then its twisted
  Jones polynomial is $-A^{-2}-A^2$ times its virtual Jones
  polynomial.
\end{thm}

If we let $M = -A^{-2}-A^2$ in the twisted Jones polynomial, we can
always divide the result by $-A^{-2}-A^2$, and then we get an
extension of the Jones polynomial to twisted links.

If the polynomial of a twisted link has an $M$ variable, then the link
is not a virtual link.  It is easy to see that if a link has an odd
number of bars on its edges, then we can factor one $M$ variable from
its polynomial.  On the right of Figure~\ref{fig:onefoil} is a onefoil
knot in a thickened Klein bottle, and its twisted Jones polynomial is:
\begin{align*}
  V_\text{Onefoil}(A,M) &= A^{-6}+(1-M^2)A^{-2}
\end{align*}
which does not have an $M$ factor.  Also its Jones polynomial is
trivial.

\begin{figure}[htbp]
  \centering
  \includegraphics[height=2in]{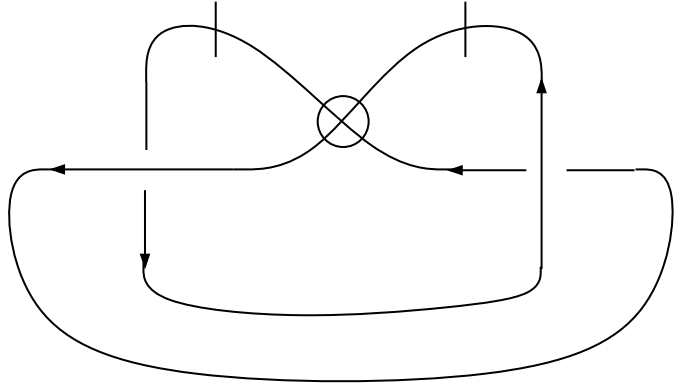}
  \caption{A non-orientable twisted knot.}
  \label{fig:T1212}
\end{figure}

However, a link may have an $M$-free twisted Jones polynomial and not
be virtual, for the knot in Figure~\ref{fig:T1212}
has twisted Jones polynomial $(-A^{-2}-A^2)(A^{-4}+A^{-6}-A^{-10})$
which is $-A^{-2}-A^2$ times its Jones polynomial.  This example is
noteworthy because while the knot is in a thickening of a projective
plane, this diagram in the projective plane is two-colorable.  In
fact, we have the following.

\begin{thm}\label{thm:link_two-colorable}
  If a twisted link has a two-colorable diagram, then its twisted
  Jones polynomial is $(-A^{-2}-A^2)$ times its Jones polynomial.
\end{thm}

\begin{figure}[htbp]
  \centering
  \includegraphics[height=2in]{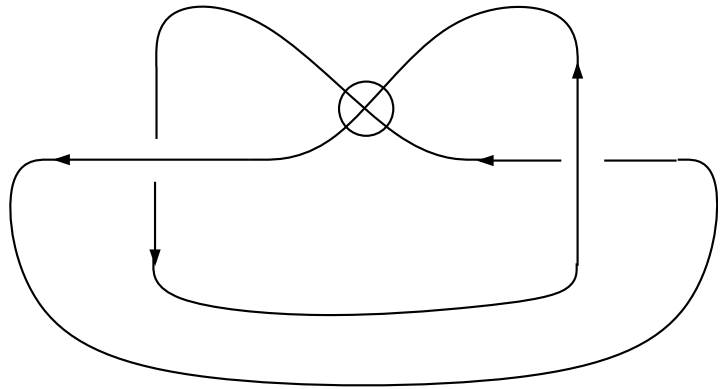}
  \caption{A virtual knot diagram that is not two-colorable.}
  \label{fig:torus1212}
\end{figure}

Again, the converse of this theorem is not true since the knot in
Figure~\ref{fig:torus1212}
is a virtual link diagram in a torus and has the same twisted Jones
polynomial as the knot in Figure~\ref{fig:T1212}, but it does not have
a two-colorable diagram because the exponents of its Jones polynomial
are not a multiple of four.  By a theorem of Naoko
Kamada~\cite{MR1914297}, if a virtual link has a two-colorable
diagram, then its Jones polynomial's exponents are multiples of four
if the link has an odd number of components, and multiples of four
plus two if the link has an even number of components.

\subsection{The Twisted Link Group}

Let $L$ be a virtual link.  Kauffman defines an invariant of virtual
links which he calls the \emph{group of the virtual link}, $\Pi
L$~\cite{MR2000i:57011}.  This invariant is calculated from any
diagram for $L$, by letting the arcs of the diagram between
undercrossings be the generators, ignoring virtual crossings, and
formally applying the Wirtinger algorithm to the diagram to create a
presentation of the group.  This group is also called the \emph{upper
  group of the virtual link}, and when the generators are arcs between
overcrossings, the group is called the \emph{lower group of the
  virtual link}.

\begin{figure}[htbp]
  \centering
  \includegraphics[height=2in]{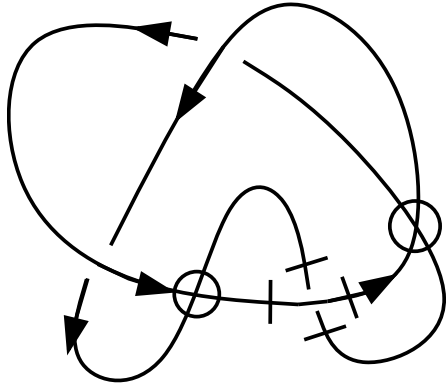}
  \caption{A trefoil under move T3.}
  \label{fig:trefoilT3}
\end{figure}

Because the group of a link ignores the bars on the edges, it fails to
be invariant with respect to move T3.  For example,
Figure~\ref{fig:trefoilT3}
shows a trefoil on which a T3 move has been performed, and its group
is trivial.  In this case, we can apply a modified Wirtinger algorithm
that defines two generators for every end of an edge in the underlying
graph, four relations at every crossing, and two relations for every
edge, and we define the \emph{twisted link group} of the link diagram
to be the group associated to the presentation.  We have the
following:

\begin{thm}\label{thm:link_twisted_group}
  The twisted link group of a link diagram is an invariant of the
  link.
\end{thm}

It is immediate from the definition of the twisted link group
that when a twisted link is virtual, its group is a free product of
its upper and lower groups.

The knots in Figure~\ref{fig:T1212} and Figure~\ref{fig:torus1212}
have identical twisted Jones polynomial, but different twisted link
groups.  The knot in Figure~\ref{fig:T1212} has a twisted link group
with (simplified) presentation:
\begin{align*}
  \twisted\Pi(\text{Twofoil}) &= \left<a,b|a^2=\inv b a^2 b\right>
\end{align*}
which has $\left<a^2\right>$ in its center.  The knot in
Figure~\ref{fig:torus1212} has a twisted link group which is a
free product of its upper and lower groups, and since these are both
isomorphic to the integers, the twisted link group has a trivial
center.

However, the twisted link group is not necessarily better than the
twisted Jones polynomial at distinguishing knots.  The knot in
Figure~\ref{fig:onefoil} and the unknot have different twisted Jones
polynomials, but identical twisted link groups.  The twisted link
group of the onefoil is a free group on two generators, as is the
twisted link group of an unknot.  But the twisted Jones polynomial of
the onefoil is:
\begin{align*}
  V_\text{Onefoil}(A,M) &= A^{-6}+(1-M^2)A^{-2}.
\end{align*}
So the twisted Jones polynomial and the twisted link group contain
different information about the link.

The topological interpretation of the twisted link group is an open
problem.  In particular, its relationship with the fundamental group
of the complement of some realization of the twisted link with
collapsed boundary is not as straightforward as in the virtual
case~\cite{bourgoin04:_fundam_group_virtual_knot}, as some examples
have shown that they are not always equal.

\subsection{Dedication}

To the memory of Jerome P.~Levine for his keen intellect, attentive
kindness and boundless generosity.

\section{Background}

\subsection{Projective Links}

\begin{figure}[htbp]
  \centering
  \raisebox{-.35in}{\includegraphics[height=.8in]{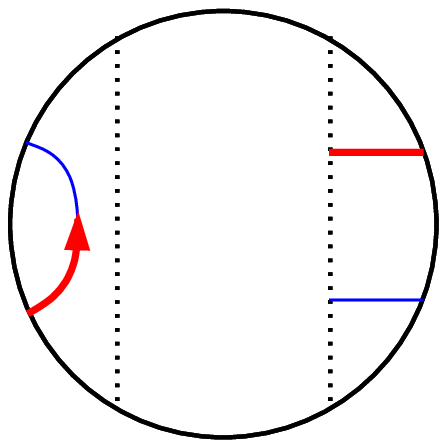}}
  $\buildrel\text{P1}\over\Longleftrightarrow$
  \raisebox{-.35in}{\includegraphics[height=.8in]{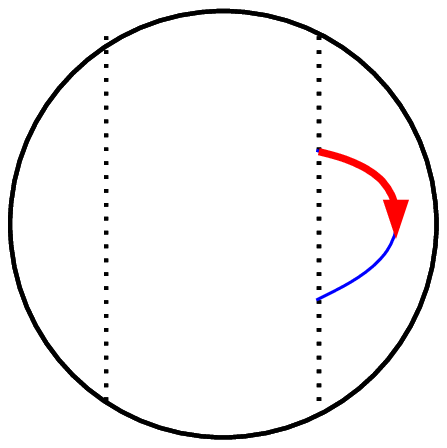}}
  \hspace{1ex}
  \raisebox{-.35in}{\includegraphics[height=.8in]{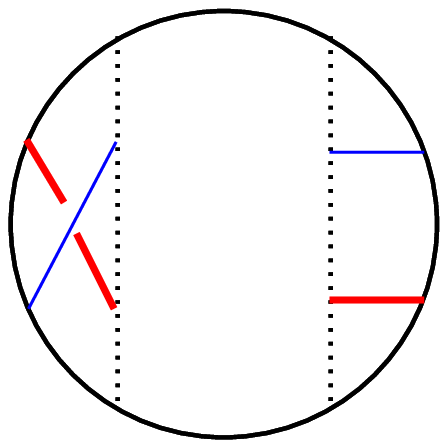}}
  $\buildrel\text{P2}\over\Longleftrightarrow$
  \raisebox{-.35in}{\includegraphics[height=.8in]{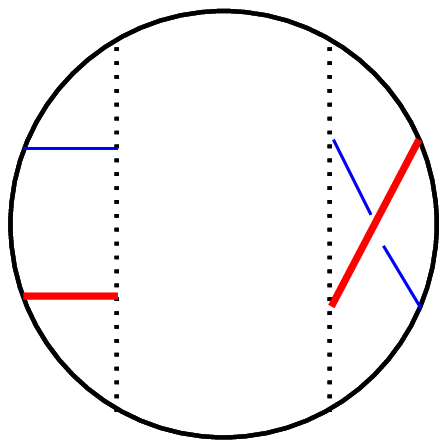}}
  \caption{Projective extension of the Reidemeister moves.}
  \label{fig:pReidemeister}
\end{figure}

In 1990, Yu.~V.~Drobotukhina introduced the study of links in real
projective space as a generalization of links in the
three-sphere~\cite{MR91i:57001}.  Drobotukhina showed how to create
diagrams of projective links in a 2-disk representation of a
projective plane, and she extended Reidemeister moves on planar
diagrams to include moves across the boundary of the disks.  On the
right of Figure~\ref{fig:twofoil}, is shown a diagram of the link in
the projective plane.  Figure~\ref{fig:pReidemeister}
shows the extension of the Reidemeister moves to deal with passing
through the boundary of the disk.  The strands involved in the moves
are shown with different thickness to emphasize the effects of passing
through the orientation-reversing boundary of the 2-disk.  If a
projective link can be isotoped to a link in the affine part of
projective space, then it corresponds to a link in the three-sphere.
In other papers, Drobotukhina classified non-trivial projective links
with diagrams with up to six crossings~\cite{MR1296890} as well as
projective Montesinos links~\cite{MR93c:57004}.

\subsection{Virtual Links}

An \emph{oriented Gauss code} is a double occurrence collection of
sequences of symbols with a direction of traversal for each sequence
and whose symbols are accompanied by writhe ($\pm$) and height ($O/U$)
marks with the restriction that both occurrences of the same symbol
have the same writhe and different heights.  Two oriented Gauss codes
over the same symbols are equivalent under permutation of the set of
symbols and rotation of the sequences.  The string $O1\p U2\p U1\p
O2+$ is an example of a one component Gauss code with two crossings
with positive writhe.  Kauffman defined a virtual link combinatorially
as follows~\cite[pg.~668]{MR2000i:57011}.

A \emph{virtual link} is an equivalence class of oriented Gauss codes
under abstractly defined Reidemeister moves for these codes shown in
Figure~\ref{fig:reidgc}.
\begin{figure}[htbp]
  \centering\small
  \begin{align*}
    &\text{R-1}:\\ &ab\ \equivalent\ aU1\p O1\p b\ \equivalent\ aU1\n O1\n b \equivalent\ aO1\p U1\p b\ \equivalent\ aO1\n U1\n b \\
    &\text{R-2}:\\ &aO1\n O2\p bU1\n U2\p c\ \equivalent\ abc\ \equivalent\ aU1\p U2\n bO1\p O2\n c \\
    &aO1\p O2\n bU2\n U1\p c\ \equivalent\ abc\ \equivalent\ aU1\n U2\p bO2\p O1\n c \\
    &\text{R-3}:\\ &aO1\p U2\p bO3\p O2\p cU1\p U3\p d\ \equivalent\ aU2\p O1\p bO2\p O3\p cU3\p U1\p d \\
    &aO1\p O2\n bO3\p U2\n cU1\p U3\p d\ \equivalent\ aO2\n O1\p bU2\n O3\p cU3\p U1\p d \\
    &aU1\n O2\n bU3\n U2\n cO1\n O3\n d\ \equivalent\ aO2\n U1\n bU2\n U3\n cO3\n O1\n d \\
    &aU1\n U2\p bU3\n O2\p cO1\n O3\n d\ \equivalent\ aU2\p U1\n
    bO2\p U3\n cO3\n O1\n d
  \end{align*}
  \caption{Reidemeister moves for oriented Gauss codes}
  \label{fig:reidgc}
\end{figure}
In the figure, the oriented Gauss code is presented left-right, the
crossing numbers are assigned to reflect the order in which they are
encountered, and the lowercase letters $a,b,c,\ldots$ represent
segments of the code in between which the fragments are embedded.

Kauffman associated every oriented Gauss code with a planar link
diagram by the usual graph theory approach of allowing drawings in the
plane where some of the crossings are deemed not to occur in the
diagram.  These crossings are called \emph{virtual} and do not have
writhe or height.  A virtual link diagram for Gauss code $O1\p U2\p
U1\p O2+$ is shown in Figure~\ref{fig:torus1212}.  In the diagram, the
virtual crossing is circled, and the writhe is obtained by the usual
right-hand rule.  Kauffman identified an extension of the Reidemeister
moves on classical link diagrams to moves on link diagrams with
virtual crossings.  Some of these are illustrated in
Figure~\ref{fig:Reidemeister} where seven of the allowed virtual moves
are shown, one for each class of Reidemeister move.  Three of the
moves involve only classical crossings, three of the moves involve
only virtual crossings, and one of the moves involves two virtual
crossings and one real crossing.  Note that the two obvious 3-crossing
moves that involve two real crossings and one virtual crossing are
forbidden because including them makes all links are equivalent to
unlinks~\cite{MR2002c:57009}.

Kauffman extended many link invariants, such as the group, quandle, to
virtual links by the simple approach of ignoring the virtual
crossings.

Abstract links over orientable surfaces were introduced in knot theory
by Naoko Kamada and Seiichi Kamada~\cite{MR2001h:57007}.  They are
regular neighborhoods of links in surfaces, and are equivalent to
virtual links which they anticipated.  However, similar presentations
for graphs embedded in surfaces have long been common in Topological
Graph Theory~\cite{MR88h:05034}.

\subsection{The Jones Polynomial}

In 1985, Vaughan Jones announced the creation of a new polynomial
invariant of links~\cite{MR86e:57006}.  This invariant encoded
different information about links than the older Alexander polynomial.
This invariant was soon generalized to the HOMFLY polynomial, which
generalized both the Jones and Alexander polynomials.  In 1987, Louis
Kauffman introduced a state-sum model for the calculation of the Jones
polynomial~\cite{MR88f:57006}.

In 1990, Yu.~V.~Drobotukhina extended the Jones polynomial to
projective links, and generalized results by L.~Kauffman and
K.~Murasugi relating the crossing and component numbers of link
diagrams to the exponents of the link's Jones
polynomial~\cite{MR91i:57001}.  Drobotukhina used the properties of
the Jones polynomial to provide a necessary condition for a projective
link to be affine, and characterized the projective links with
alternating diagrams that are affine.  Recently, M.~Mroczkowki
considered the problem of creating descending diagrams for projective
links~\cite{MR1999636}, and used them to extend the HOMFLY and
Kauffman polynomials to projective links~\cite{MR2128052}.

In 1996, Louis Kauffman extended the Jones polynomial to virtual links
by the simple approach of ignoring the virtual
crossings~\cite{MR2000i:57011}.  For example, the bracket polynomial
of the virtual link in Figure~\ref{fig:torus1212} is $A^{2}+1+A^{-4}$
showing that it is not trivial.  Soon afterward, the Alexander,
HOMFLY, and Kauffman polynomials were also extended to virtual links.

\subsection{The Group of a Virtual Link}

In~\cite{MR2000i:57011}, Louis Kauffman introduced the group of a
virtual link as the group defined by the formal Wirtinger presentation
obtained from any link diagram by ignoring any virtual crossings.
Unlike the classical group of a link, the group of a virtual link can
have deficiency zero.  More importantly, it may not be residually
finite, so it may not be the group of any
three-manifold~\cite{MR2002m:57011}, and in particular, it may not be
the fundamental group of the complement of the link in any of its
stabilized embeddings.  Furthermore, the group obtained by defining
the generators to be strands between over-crossings may be different
from the standard one where generators are strands between
under-crossings.

In~\cite{MR2001h:57007}, Naoko Kamada and Seiichi Kamada showed that
the group of a virtual knot is the group of a three-complex obtained
by collapsing one of the boundary components of the complement of the
link in the thickened surface.

\section{Definitions}

The \emph{Euler genus} of a surface is $2-\chi$ where $\chi$ is the
surface's Euler characteristic.  A closed surface of Euler genus $g$
will be noted $\Sigma_g$, a closed orientable surface of orientable
genus $g$ will be noted $S_g$, and a closed non-orientable surface of
Euler genus $g$ will be noted $N_g$.

A \emph{move} of a manifold $M$ is a homeomorphism of $M$ supported by
a ball, keeping the boundary of the ball fixed.  A \emph{standard
  linear move} is a homeomorphism of a standard simplex keeping the
boundary fixed, mapping the barycenter to another interior point, and
joining linearly.  A \emph{linear move} of $M$ is a move $h$ supported
by a ball $B$ for which there exists a homeomorphism $k:B\to M$ such
that $kh\inv k$ is a standard linear move.  Two embeddings $f,g:N\to
M$ are \emph{isotopic by linear moves} if there exists a finite
sequence $h_1,h_2,\ldots,h_n$ of linear moves such that $h_1h_2\cdots
h_nf = g$.

Define the standard interval to be $I = [-1,1]$.  An \emph{orientable
  thickening} of a closed surface $\Sigma$ is the orientation
$I$-bundle $\Sigma\ttimes I$ over the surface.  An \emph{oriented
  thickening} is an orientable thickening with a given orientation.  A
\emph{link in an oriented thickening} is an embedding of a disjoint
collection of oriented circles in the interior of an oriented
thickening.  Two links in an oriented thickening are \emph{equivalent}
if there is an orientation-preserving homeomorphism of the oriented
thickening that is isotopic to the identity and takes one link to the
other while preserving link orientation.

A \emph{band} is either an annulus or a M\"obius band.  A
\emph{vertical band} in an oriented thickening is the fiber over a
simple closed path in the zero-section.  A \emph{vertical annulus} is
a vertical band over an orientation-preserving path, and a
\emph{vertical M\"obius band} is a vertical band over an
orientation-reversing path.  A \emph{topologically vertical band} is a
properly embedded band isotopic to a vertical band.

A \emph{destabilization} of a link in an oriented thickening consists
of cutting the thickening along a vertical band disjoint from the link
and capping the resulting boundary components with thickened disks.
The result of a destabilization is a link in an oriented thickening
\emph{descendant} from the original link.  When the band is a vertical
M\"obius band, it has a non-trivial normal bundle in the zero-section
so its normal bundle in the oriented thickening is trivial.  Then, the
boundary component that remains after cutting along it is an annulus
that is capped with a single thickened disk.  A \emph{stabilization}
is the reverse of a destabilization.  Two destabilizations are
\emph{descent equivalent} if they have equivalent descendants.

\begin{figure}[htbp]
  \centering
  \includegraphics[height=1.5in]{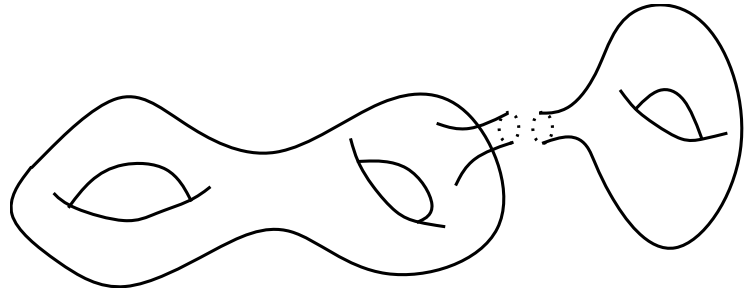}
  \caption{A (de)stabilization with a thickened handle.}
  \label{fig:stablehandle}
\end{figure}

\begin{figure}[htbp]
  \centering
  \includegraphics[height=1.5in]{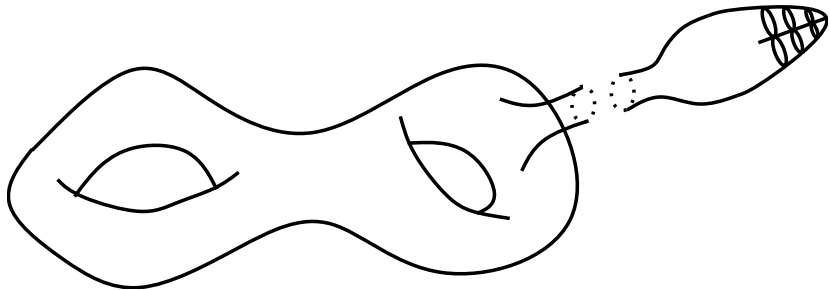}
  \caption{A (de)stabilization with a thickened crosscap.}
  \label{fig:stablecrosscap}
\end{figure}

Figure~\ref{fig:stablehandle}
shows the (de)stabilization of a thickened surface (removing or)
adding a thickened handle.  Figure~\ref{fig:stablecrosscap}
shows the (de)stabilization of a thickened surface (removing or)
adding a thickened crosscap.

A \emph{stable link in an oriented thickening} is an equivalence class
of links in oriented thickenings under (de)stabilizations.  Unless
explicitly said otherwise, a mention of links in oriented thickenings
will always imply equivalence up to stabilization.

\begin{figure}[htbp]
  \centering
  \includegraphics[height=1.5in]{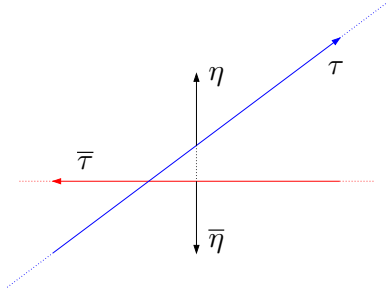}
  \caption{Resolving a double point of a surface projection.}
  \label{fig:Xpos}
\end{figure}

A \emph{surface projection} of a link in an oriented thickening is a
regular projection of the link by the bundle projection on the
zero-section of the oriented thickening.  The components of a surface
projection are oriented by the components of the link.  Let $\tau,
\other\tau$ be unit vectors at a double point each tangent to
different parts of the curve through the double point, and that agree
with the orientation of the part to which they are tangent.  A
\emph{resolution} of a double point of a surface projection is an
assignment of unit vectors $\eta, \other\eta$ from the normal bundle
to the surface at the double point such that $\eta = -\other\eta$.
Figure~\ref{fig:Xpos}
shows the resolution of a double point.  Then, the \emph{writhe}
$\writhe$ at a double point of a generic immersion is $\writhe =
(\tau\times\other\tau)\cdot\eta=(\other\tau\times\tau)\cdot\other\eta$.

\begin{figure}[htbp]
  \centering
  \raisebox{-.75in}{\includegraphics[height=1.5in]{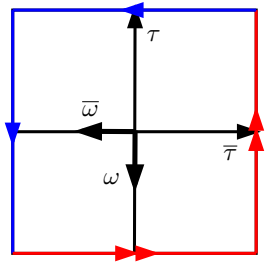}}
  \caption{A link surface diagram in a Klein bottle.}
  \label{fig:Sonefoil}
\end{figure}

We now consider generic immersions of oriented curves in compact
surfaces as in~\cite{bourgoin03:_class_immer_curves}.  Let $\tau,
\other\tau$ be unit vectors at a double point each tangent to
different parts of the curve through the double point, and that agree
with the orientation of the part to which they are tangent.  A
\emph{separation} of a double point of a generic immersion is an
assignment of unit vectors $\omega, \other\omega$ at the double point
such that $\omega = \pm\tau, \other\omega = \pm\other\tau$ and such
that $\tau\cdot\omega = \other\tau\cdot\other\omega$.  Then, the
\emph{writhe} $\writhe$ at a double point of a generic immersion is
$\writhe = \tau\cdot\omega$.  A \emph{link surface diagram} is a generic
immersion with a separation of each double point.
Figure~\ref{fig:Sonefoil}
shows a link surface diagram for the onefoil in a Klein bottle with
the short fat arrows indicating the separation, in which case, the
crossing has -1~writhe.  As will be shown later separations of double
points define resolutions of the double points once the surface is
given an oriented thickening.  And unlike the usual over-under knot
diagram crossings, they do not need to change if the crossing is
isotoped along an orientation-reversing loop in the surface.  In
drawings of link diagrams in surfaces, the separation can be indicated
by the classical over--under convention, so long as an isotopy of a
crossing through an orientation-reversing path switches the crossing.

\begin{figure}[htbp]
  \centering
  \raisebox{-.35in}{\includegraphics[height=.8in]{1vert}}    
  $\buildrel\text{R1}\over\Longleftrightarrow$
  \raisebox{-.35in}{\includegraphics[height=.8in]{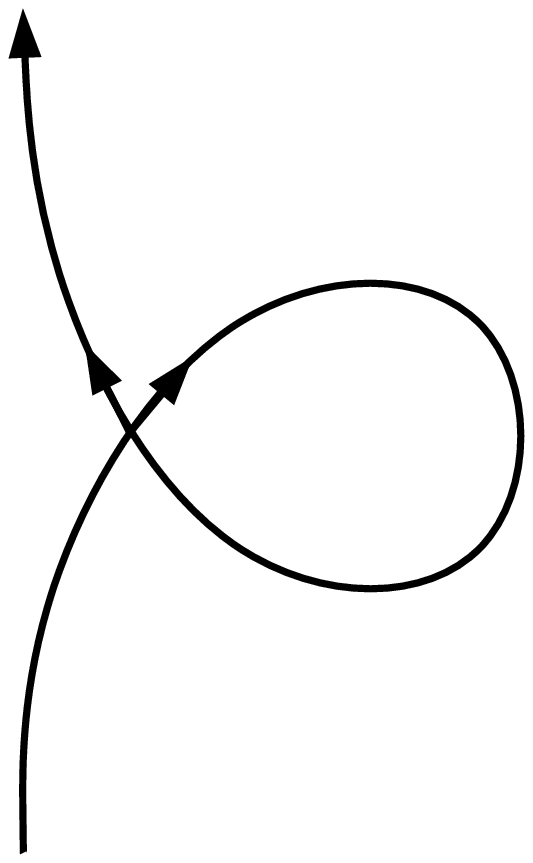}}
  \hspace{1ex}
  \raisebox{-.35in}{\includegraphics[height=.8in]{2vert}}
  $\buildrel\text{R2}\over\Longleftrightarrow$
  \raisebox{-.35in}{\includegraphics[height=.8in]{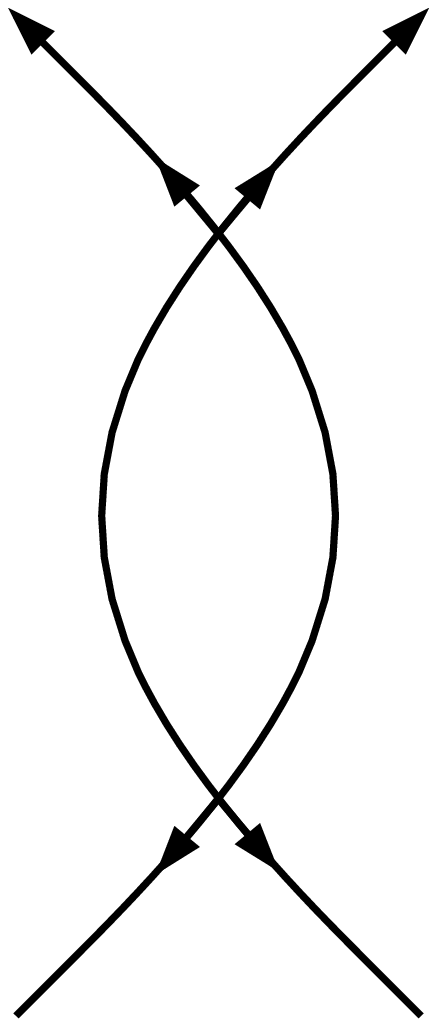}}
  \hspace{1ex}
  \raisebox{-.35in}{\includegraphics[height=.8in]{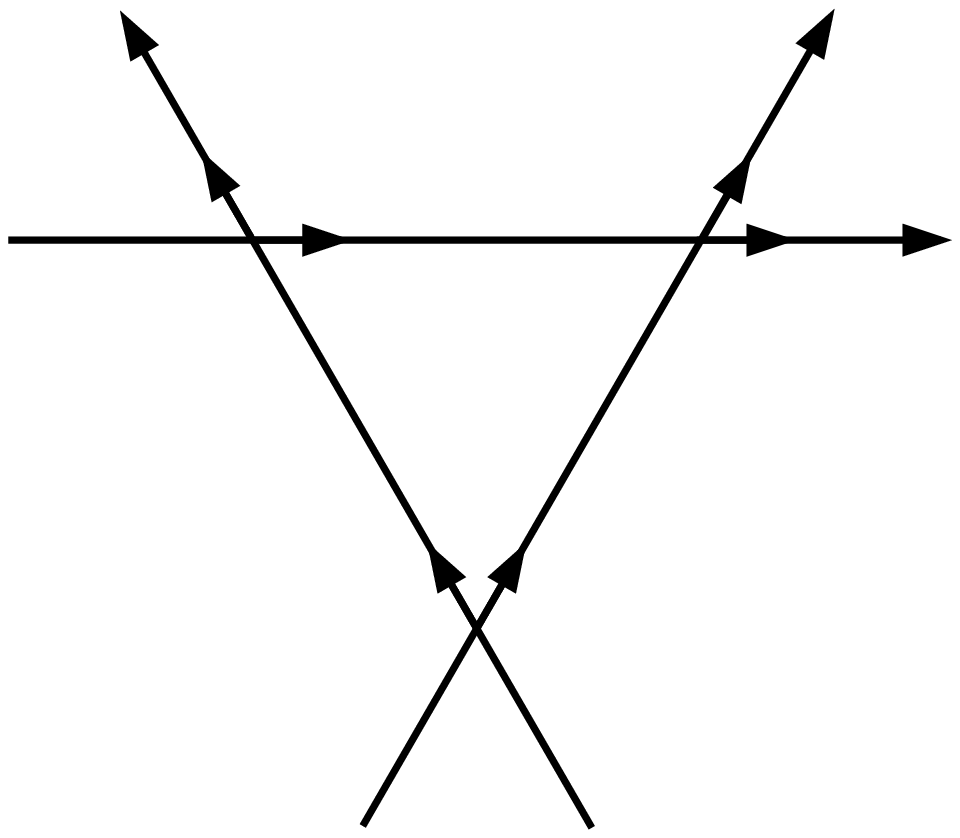}}
  $\buildrel\text{R3}\over\Longleftrightarrow$
  \raisebox{-.35in}{\includegraphics[height=.8in]{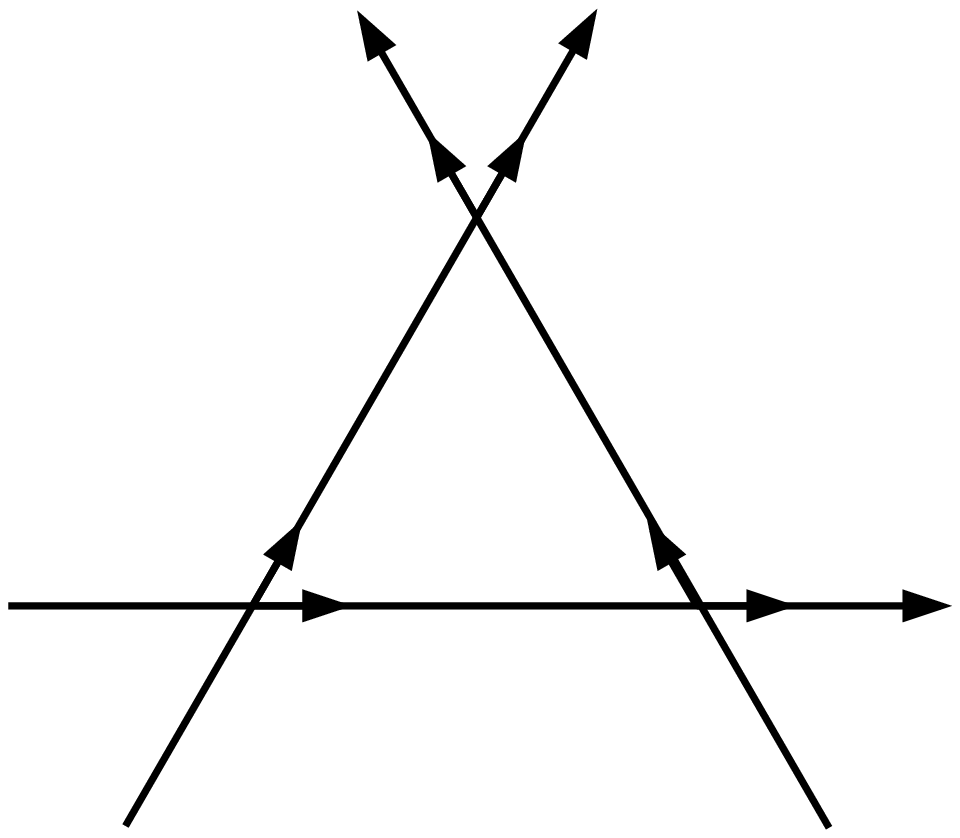}} \\[1ex]
  \caption{Reidemeister moves for link surface diagrams.}
  \label{fig:lsd_Reidemeister}
\end{figure}

A \emph{Reidemeister move} on a link surface diagram is one of the
three moves illustrated in Figure~\ref{fig:lsd_Reidemeister}.
Two link surface diagrams in homeomorphic surfaces are
\emph{Reidemeister equivalent} if there are sequences of ambient
isotopies of the surfaces and Reidemeister moves in disk neighborhoods
such that there exists a homeomorphism of the surfaces that carries
one diagram onto the other one while preserving the diagrams'
orientations.

\begin{figure}[htbp]
  \centering
  \raisebox{-.25in}{\includegraphics[width=1in]{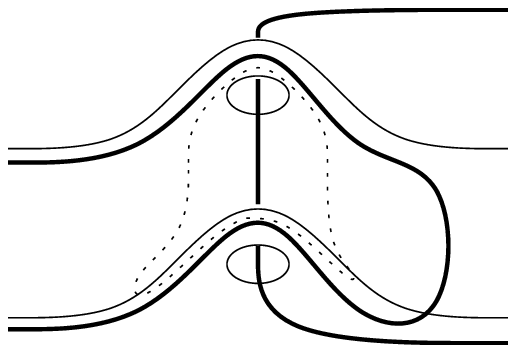}}
  $\Leftrightarrow$
  \raisebox{-.1in}{\includegraphics[width=1in]{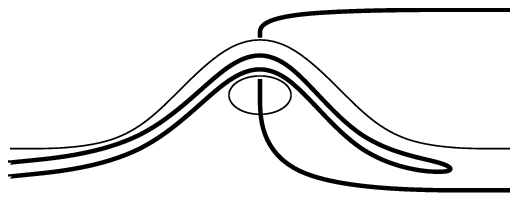}}
  $\Leftrightarrow$
  \raisebox{-.1in}{\includegraphics[width=1in]{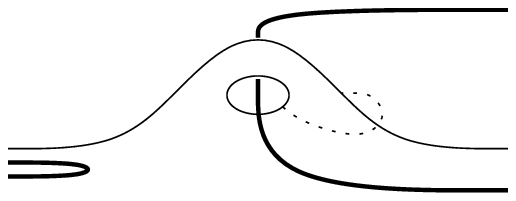}}
  $\Leftrightarrow$
  \raisebox{-.1in}{\includegraphics[width=1in]{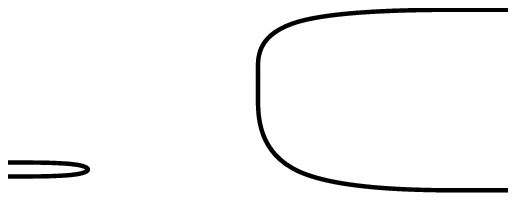}}
  \caption{Two surface (de)stabilizations}
  \label{fig:v2}
\end{figure}

A \emph{destabilization} of a link surface diagram consists of cutting
the surface along a simple closed path disjoint from the link and
capping the resulting boundary components with disks.  The result of a
destabilization is a link surface diagram \emph{descendant} from the
original link.  When the path is orientation-reversing, the boundary
component that remains after cutting along it is a circle that is
capped with a single disk.  A \emph{stabilization} is the reverse of a
destabilization.  Two destabilizations are \emph{descent equivalent}
if they have Reidemeister-equivalent descendants.  A simple closed
path is \emph{essential} if it does not bound a disk in the complement
of a link surface diagram.  Figure~\ref{fig:v2}
shows a sequence of two destabilizations along dotted paths that are
essential and remove two handles.  This is equivalent to performing a
V2-move on a virtual link diagram.  A \emph{stable link surface
  diagram} is an equivalence class of link surface diagrams under
(de)stabilizations of the surface that do not intersect the link.  Two
stable link surface diagrams are \emph{stably Reidemeister equivalent}
if there is a sequence of (de)stabilizations of the surfaces,
homeomorphisms of the (de)stabilized surfaces, ambient isotopies of
the surfaces, and Reidemeister moves in disk neighborhoods that
carries one diagram onto the other one.  Unless explicitly said
otherwise, ``link surface diagram'' will always mean stable link
surface diagram.

An \emph{intersigned link code} is a double occurrence collection of
sequences of symbols with a direction of traversal for each sequence
and whose symbols are accompanied by writhe ($\pm$) and separated by
signs~($\pm$) with the restriction that both occurrences of the same
symbol have the same writhe.  They are in fact intersigned Gauss codes
with writhe marks attached to the symbols, and intersigned Gauss codes
classify immersed curves in closed compact
surfaces~\cite{bourgoin03:_class_immer_curves}.  A \emph{Reidemeister
  move} on an intersigned link code is an abstract Reidemeister moves
described in Figure~\ref{fig:reidilc}.
\begin{figure}[htbp]
  \centering
  \begin{align*}
    &\text{R-1}: &&ab~\equivalent~a 1^\wm\p1^\wm b \\
    &\text{R-2}: &&a 1^\wm\p2^{-\wm} b 1^\wm\p2^{-\wm}
    c~\equivalent~abc~\equivalent~a 1^\wm\p2^{-\wm} b 2^{-\wm}\p1^\wm c
    \\
    &\text{R-3}: &&a 1^\wm\p2^\wm b 3^\wm\n2^\wm c 1^\wm\n3^\wm
    d~\equivalent~a 1^\wm\p2^\wm b 1^\wm\n3^\wm c 3^\wm\n2^\wm d \\
    &&&a 1^{-\wm}\p2^\wm b 3^\wm\n2^\wm c 1^{-\wm}\n3^\wm d~\equivalent~a
    1^\wm\p2^{-\wm} b 1^\wm\n3^\wm c 3^\wm\n2^{-\wm} d
  \end{align*}
  \caption{Reidemeister moves for intersigned link codes}
  \label{fig:reidilc}
\end{figure}
In the figure, the intersigned link code is presented left-right, the
crossing numbers are assigned to reflect the order in which they are
encountered, the writhe mark is given as an exponent $\wm = \pm$ such
that $-\wm = \mp$, and the lowercase letters $a,b,c,\ldots$ represent
segments of the code in between which the fragments are embedded.  Two
intersigned link codes are \emph{Reidemeister equivalent} if there is
a sequence of Reidemeister moves that takes one code to the other.

\begin{figure}[htbp]
  \centering
  \raisebox{-.5in}{\includegraphics[height=1.5in,angle=270,origin=c]{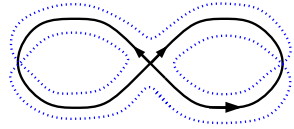}}
  \caption{An abstract knot for the onefoil.}
  \label{fig:onefoilA}
\end{figure}

\begin{figure}[htbp]
  \centering
  \raisebox{-.75in}{\includegraphics[height=1.5in]{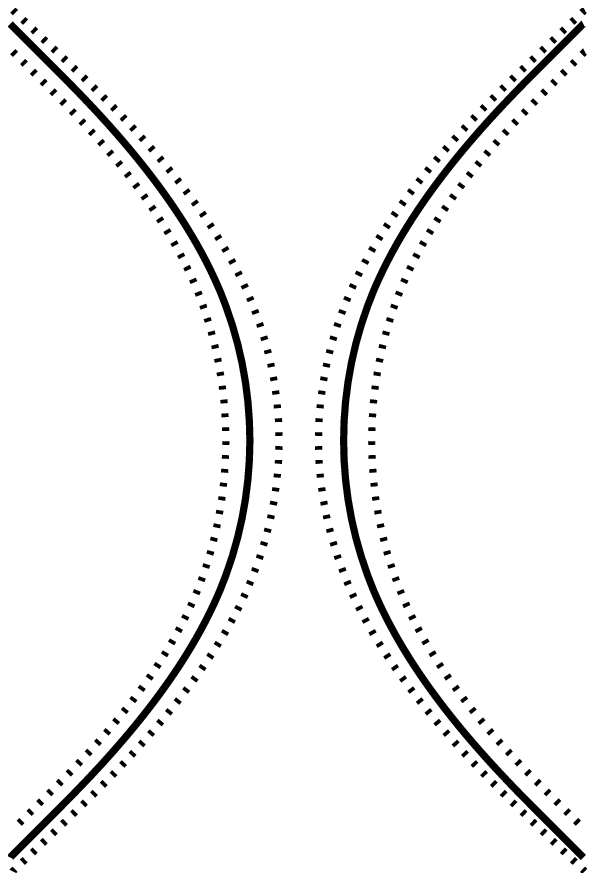}}
  $\Leftrightarrow$
  \raisebox{-.75in}{\includegraphics[height=1.5in]{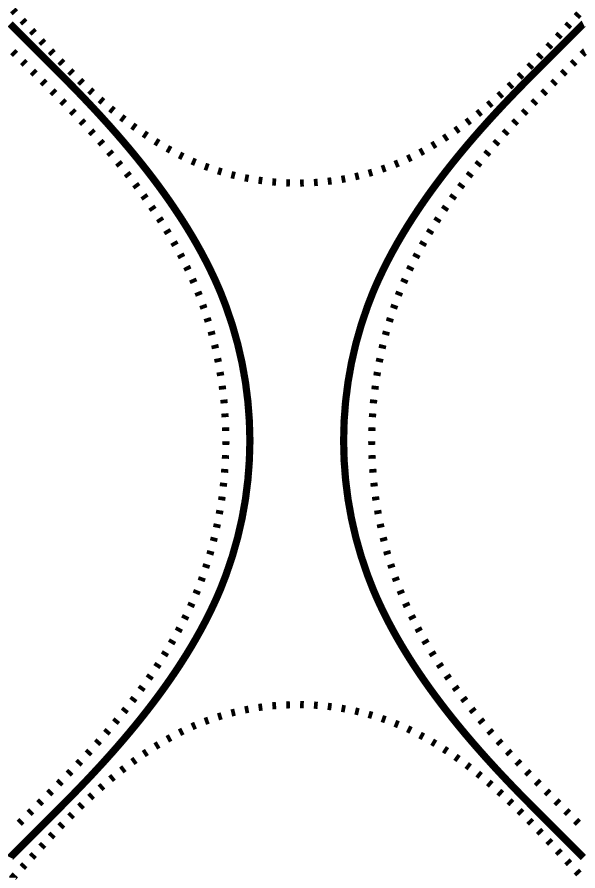}}
  $\Leftrightarrow$
  \raisebox{-.75in}{\includegraphics[height=1.5in]{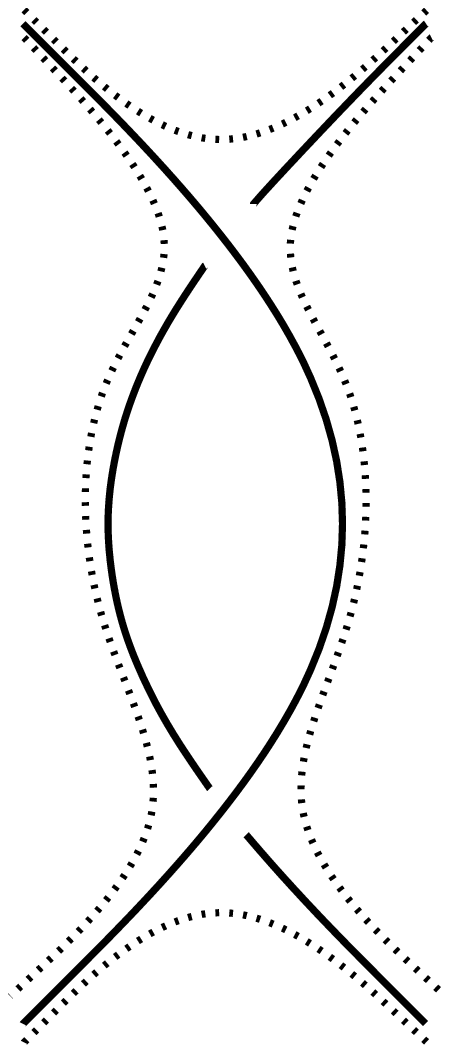}}
  $\Leftrightarrow$
  \raisebox{-.75in}{\includegraphics[height=1.5in]{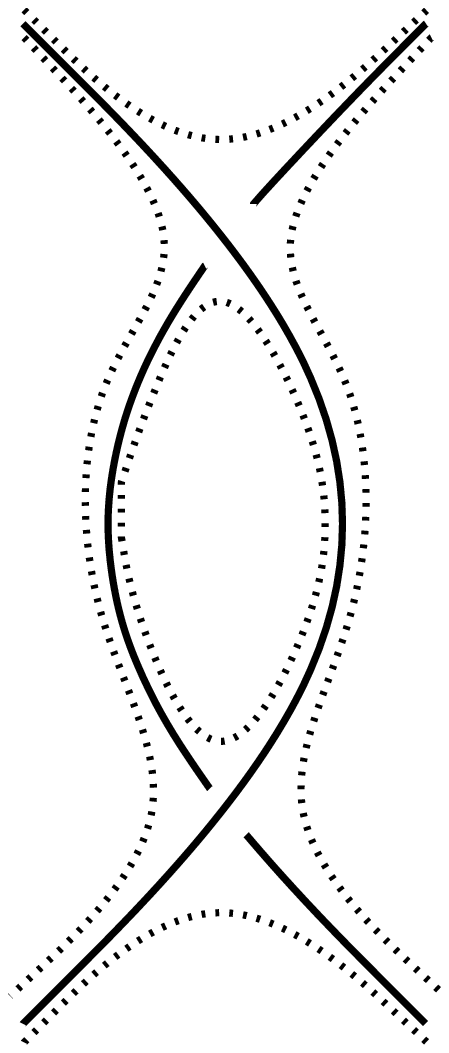}}
  \caption{An R2 move on an abstract link}
  \label{fig:abstractR2}
\end{figure}

A \emph{ribbon graph} is a finite collection of oriented disks with a
finite collection of disjoint arcs in their boundaries and (possibly
orientation-reversing) homeomorphisms that identify pairs of the
boundary arcs such that each boundary arc is identified with exactly
one other different boundary arc.  Then, the boundary of a ribbon
graph is a collection of disjoint orientation-preserving closed paths
called the \emph{faces} of the ribbon graph.  An \emph{abstract link}
is a link surface diagram in the image of an embedding of a ribbon
graph in a closed surface such that the link is a deformation retract
of the image of the ribbon graph.  This definition of an abstract link
diagram differs from that of Naoko Kamada and Seiichi
Kamada~\cite{MR2001h:57007} in that in our case, the ribbon graph need
not be an orientable surface.  Figure~\ref{fig:onefoilA}
shows an abstract link for the onefoil knot in a Klein bottle minus a
disk.  A \emph{Reidemeister move} on an abstract link consists of
first identifying intervals in the boundary of the ribbon graph to
form a disk neighborhood of a collection of edges of the link,
performing the Reidemeister move in that disk, and then discarding all
but a regular neighborhood of the link.  Figure~\ref{fig:abstractR2}
shows an R2 move on an abstract link.  Two abstract links are
\emph{Reidemeister equivalent} if there is a sequence of Reidemeister
moves that will make them homeomorphic while preserving link
orientation.

\begin{figure}[htbp]
  \centering
  \includegraphics[height=1.5in]{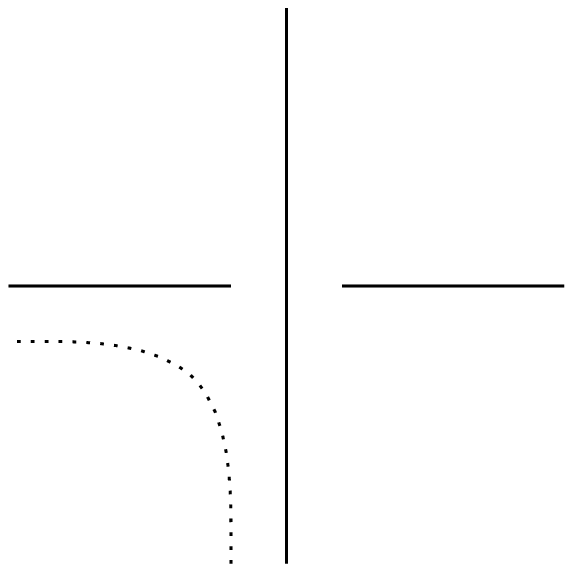}
  \hspace{1ex}
  \includegraphics[height=1.5in]{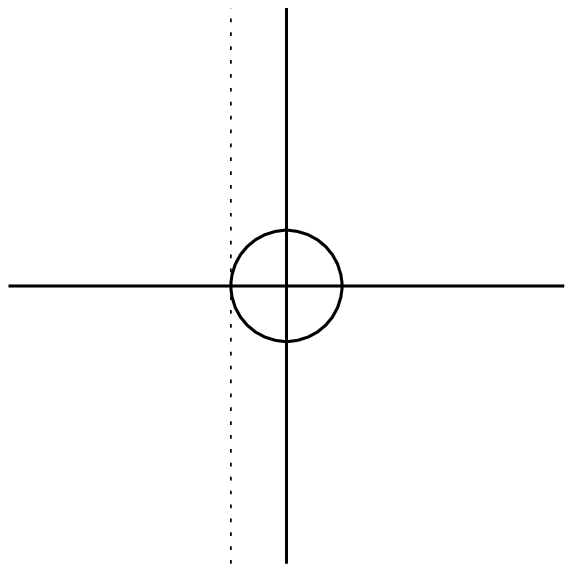}
  \hspace{1ex}
  \includegraphics[height=1.5in]{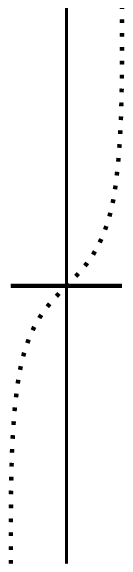}
  \hspace{1ex}
  \caption{Faces and crossings, virtual crossings, and bars.}
  \label{fig:facesand}
\end{figure}

A \emph{twisted link diagram} is a planar immersion of a disjoint
union of oriented circles with real and virtual crossings and bar
marks on edges.  Two twisted link diagrams are \emph{Reidemeister
  equivalent} if there is a sequence of ambient isotopies of the plane
and extended Reidemeister moves that carries one diagram onto the
other one.  The ten extended Reidemeister moves are diagrammed in
Figure~\ref{fig:Reidemeister}.  The \emph{faces} of a twisted link
diagram are closed curves that run along the immersed curve and have
the relationship with the crossings, virtual crossings, and bars as
show in Figure~\ref{fig:facesand}.
At a crossing, a face turns so as to avoid crossing the link diagram.
At a virtual crossing, a face goes through the virtual crossing.  At a
bar, a face crosses to the other side of the link diagram.  A twisted
link diagram is \emph{two-colorable} if its faces can be assigned one
of two colors such that the arcs of the link diagram between two
crossings always separate faces of one color from those of the other.

\section{Unique Destabilization}

According to Greg Kuperberg~\cite{MR1997331}, links in oriented
thickenings of orientable closed surfaces have unique irreducible
representatives.  This is proven by showing that all pairs of
destabilizations along essential vertical annuli are descent
equivalent to destabilizations along annuli that do not intersect, and
these are descent equivalent.
Theorem~\ref{thm:unique_destabilization} is the equivalent statement
for links in oriented thickenings over closed surfaces, and our proof
will follow Kuperberg's strategy.

\begin{proof}[Proof of Theorem~\ref{thm:unique_destabilization}]
  We expand the range of surfaces along which destabilizations can
  occur to include spheres and proper disks that separate some
  components of the link from the other components, and bands that
  separate some link-free part of the genus from the surface and so
  that we may discard that part.  The results of such destabilizations
  can be achieved by destabilizations along vertical bands.  In the
  first case, a disk or sphere may be altered to be a band, and in the
  second case, the destabilization may be performed along bands that
  progressively reduce genus but do not increase the number of
  components of the thickened surface.  So we define an
  \emph{admissible surface} to be a vertical band, sphere, or proper
  disk.  And an admissible surface is \emph{essential} if it does not
  bound a ball in the complement of the link.  Then, every essential
  admissible surface can be used for destabilization.

  Let $L\subset\Sigma_g\ttimes I = M$ be an $n$-component link in the
  oriented thickening of a surface $\Sigma_g$ of Euler genus $g$ and
  having $c$ components.  We require that every component of $M$
  contain a part of the link, so $n\ge c$.
  
  We induce on the complexity of the intersection between two
  destabilization surfaces.  Suppose that $L$ is a link with minimal
  $g+2(n-c)$ such that it has at least two descent-inequivalent
  destabilization surfaces, $B_1, B_2$, and that these intersect in
  the fewest curves amongst such pairs.  This intersection is not
  empty because destabilizations along disjoint surfaces are descent
  equivalent.  Indeed, if the surfaces are isotopic, then
  destabilization along one or the other produces the same manifold,
  and otherwise, we may destabilize first along either surface and
  then along the other to obtain the same manifold.  We consider the
  surfaces to be in general position, so their intersection is an
  embedded one-manifold, and its components are, up to homotopy, one
  of the types shown in Figure~\ref{fig:component_types}.
  \begin{figure}[htbp]
    \centering
    \parbox{.8in}{\centering\includegraphics[height=.8in]{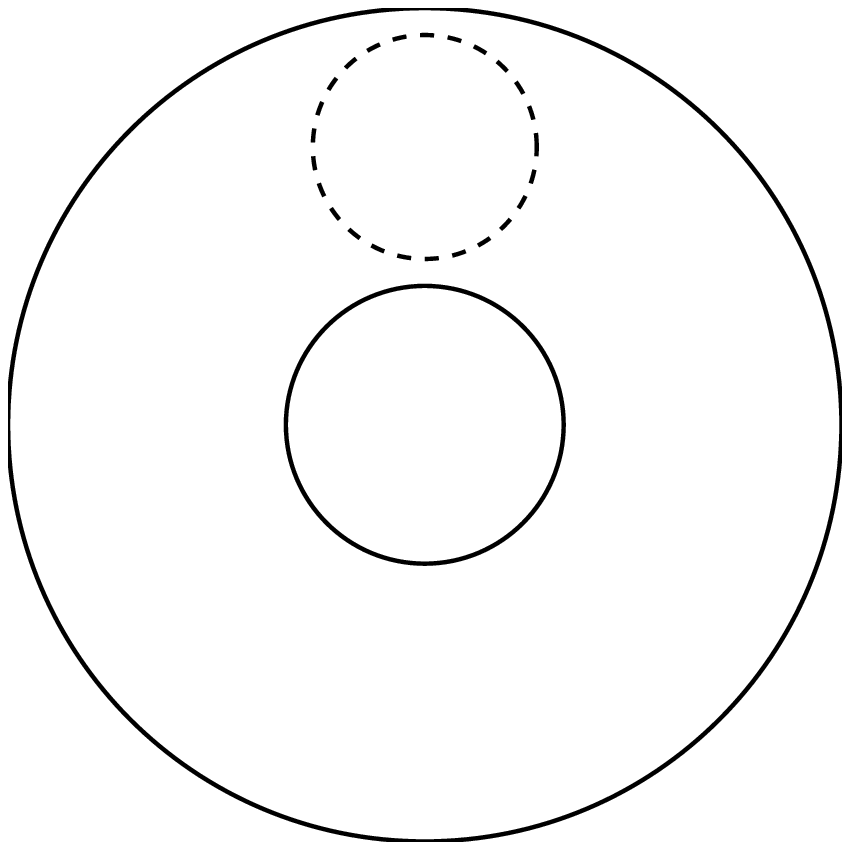}\\a1}\hspace{1em}
    \parbox{.8in}{\centering\includegraphics[height=.8in]{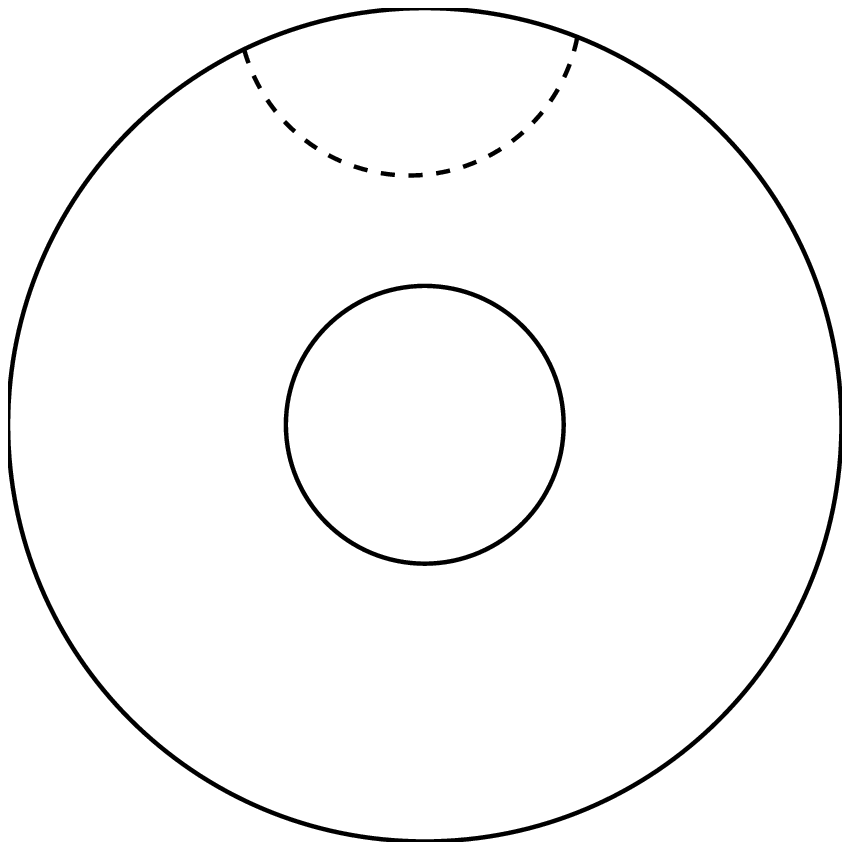}\\a2}\hspace{1em}
    \parbox{.8in}{\centering\includegraphics[height=.8in]{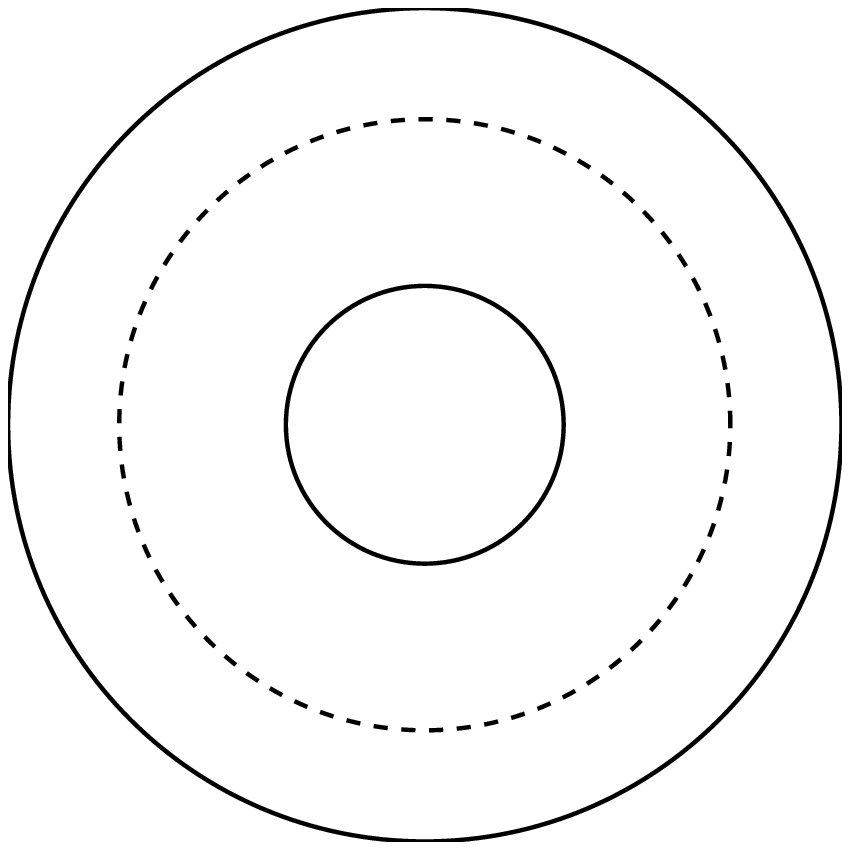}\\a3}\hspace{1em}
    \parbox{.8in}{\centering\includegraphics[height=.8in]{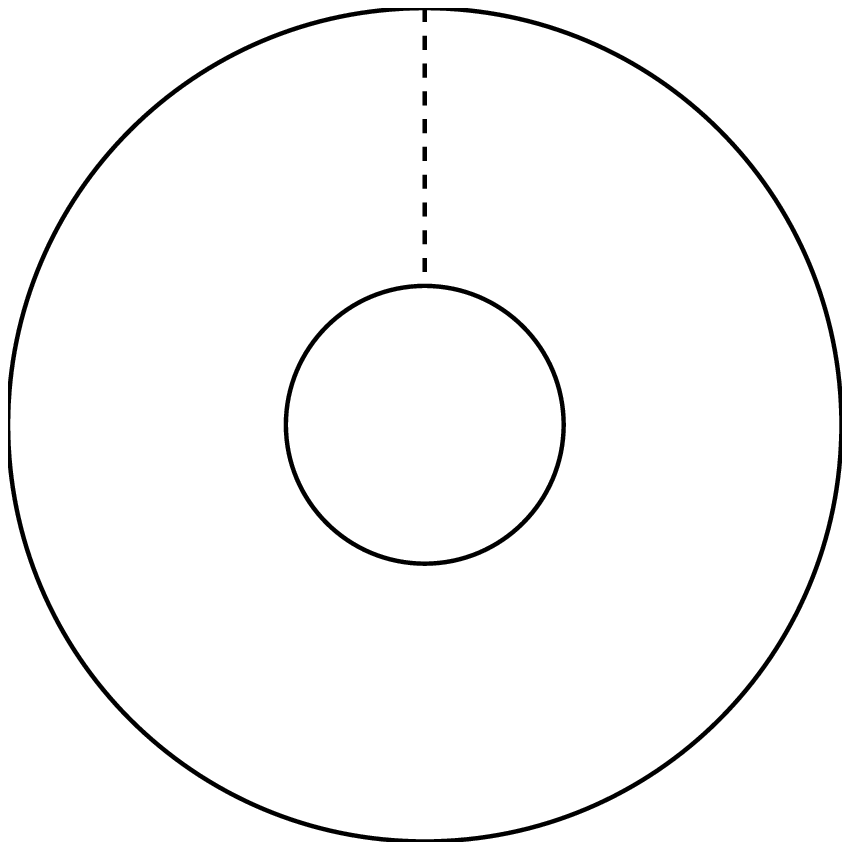}\\a4}\\[2ex]
    \parbox{.8in}{\centering\includegraphics[height=.8in]{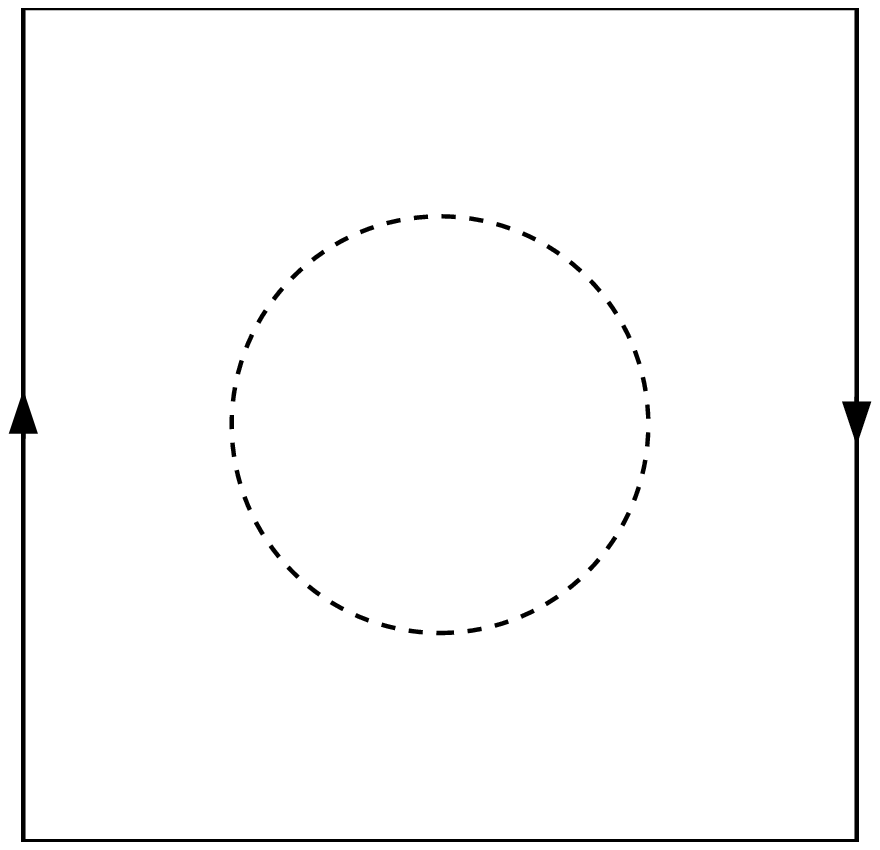}\\m1}\hspace{1em}
    \parbox{.8in}{\centering\includegraphics[height=.8in]{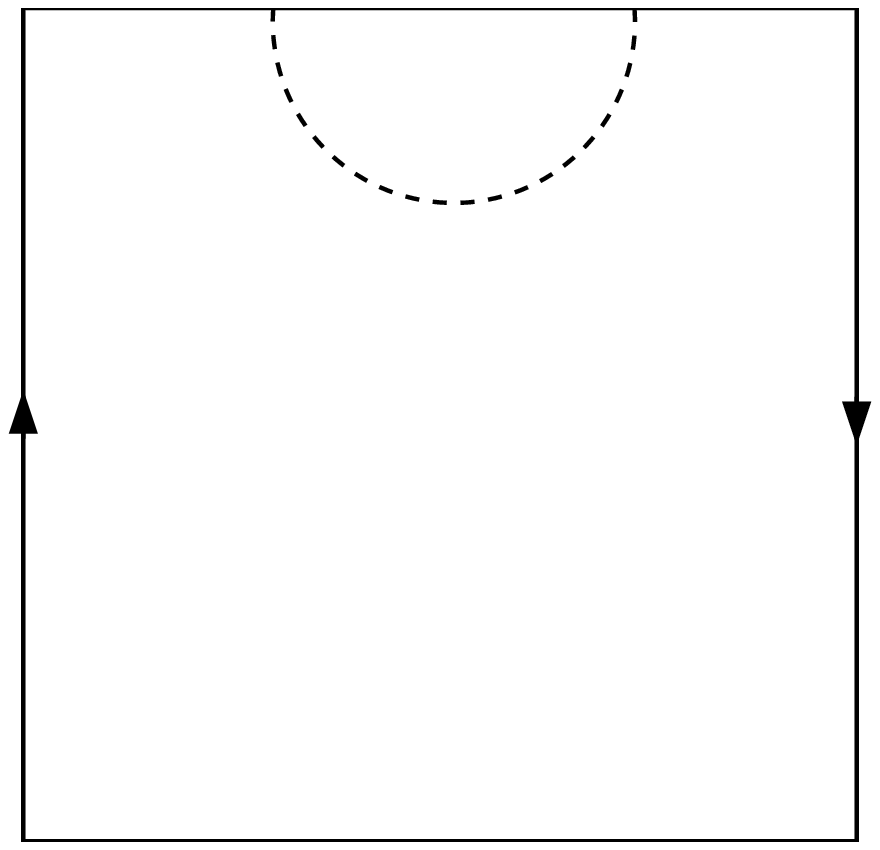}\\m2}\hspace{1em}
    \parbox{.8in}{\centering\includegraphics[height=.8in]{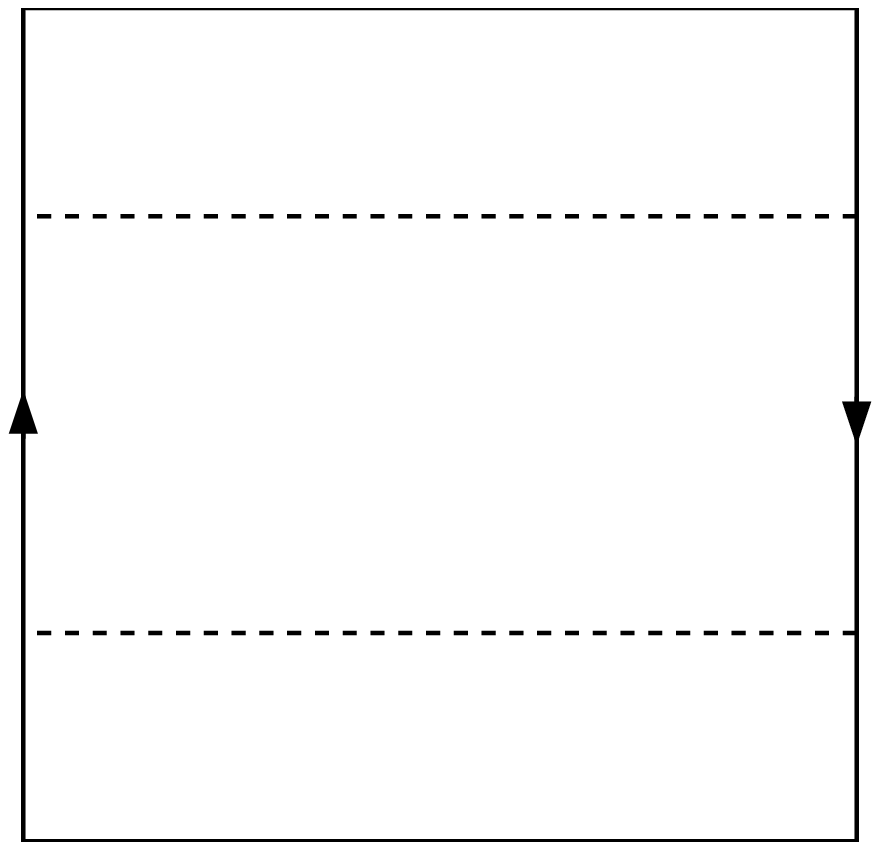}\\m3}\hspace{1em}
    \parbox{.8in}{\centering\includegraphics[height=.8in]{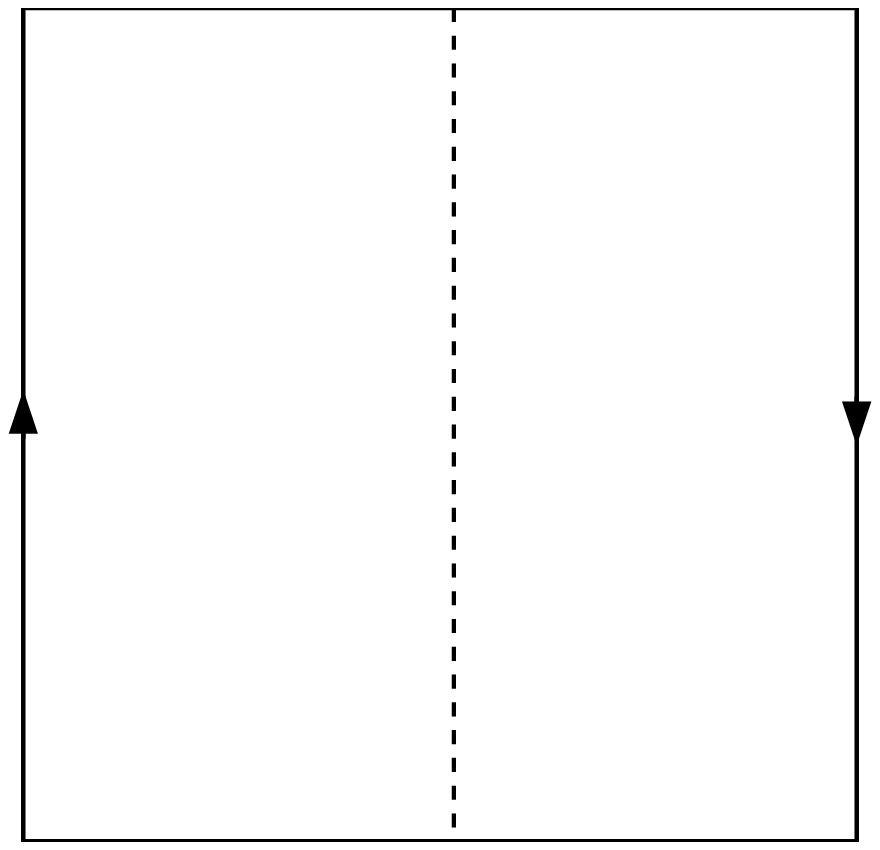}\\m4}\hspace{1em}
    \parbox{.8in}{\centering\includegraphics[height=.8in]{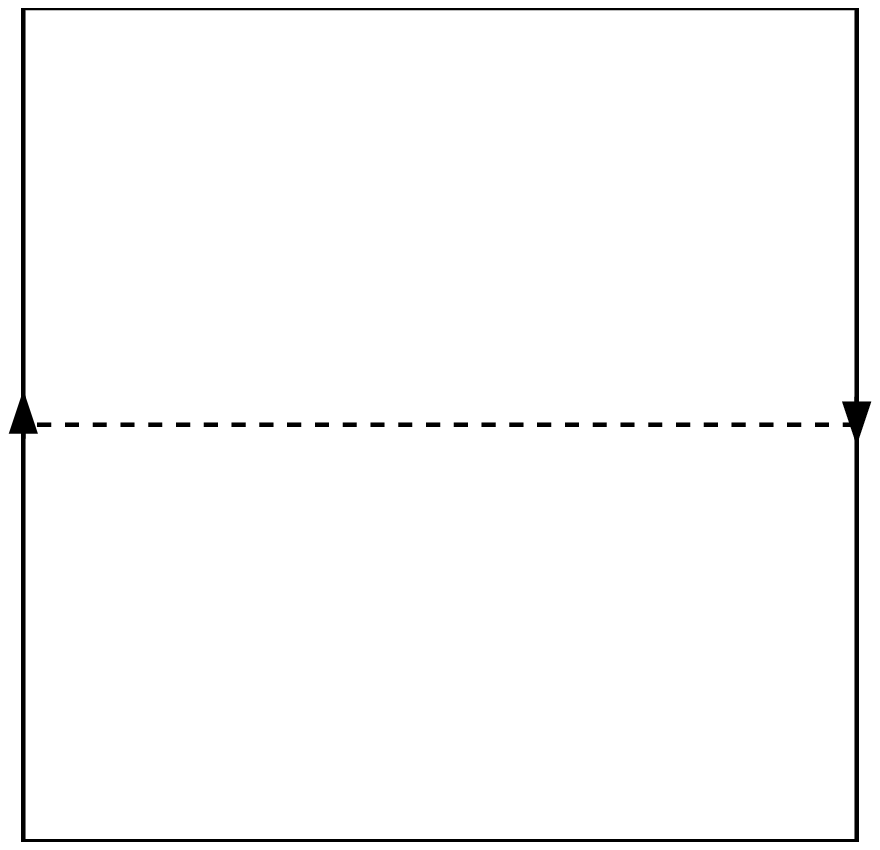}\\m5}
    \caption{Types of components of the intersection of two surfaces.}
    \label{fig:component_types}
  \end{figure}
  We deal with each type in turn.  We show that using an intersection
  component, we can find an essential destabilization surface $B'$
  whose intersection with $B_1$ and $B_2$ has fewer curves with either
  than they have with each other, so it is descent equivalent to both,
  thus obtaining a contradiction.  When $B'$ is to be a band
  constructed from parts of $B_1$ and $B_2$, we will obtain that it is
  topologically vertical because it is homotopic to the fiber over a
  closed path in the zero section made from parts of the closed paths
  of $B_1$ and $B_2$.
  
  \emph{Cases a1, a2, m1, and m2}: The curve is an arc or circle in
  $B_1$.  Then we can find an \emph{innermost} curve that bounds a
  disk that does not contain any other component of the intersection.
  (In cases a2 and m2, the disk is formed with part of the boundary of
  $B_1$.)  We may compress $B_2$ along this disk to get surfaces
  $B'_2, B''_2$ at least one of which, say $B'_2$ is an admissible
  essential surface.  In Figure~\ref{fig:compression}
  \begin{figure}[htbp]
    \centering
    \raisebox{-.75in}{\includegraphics[height=1.5in]{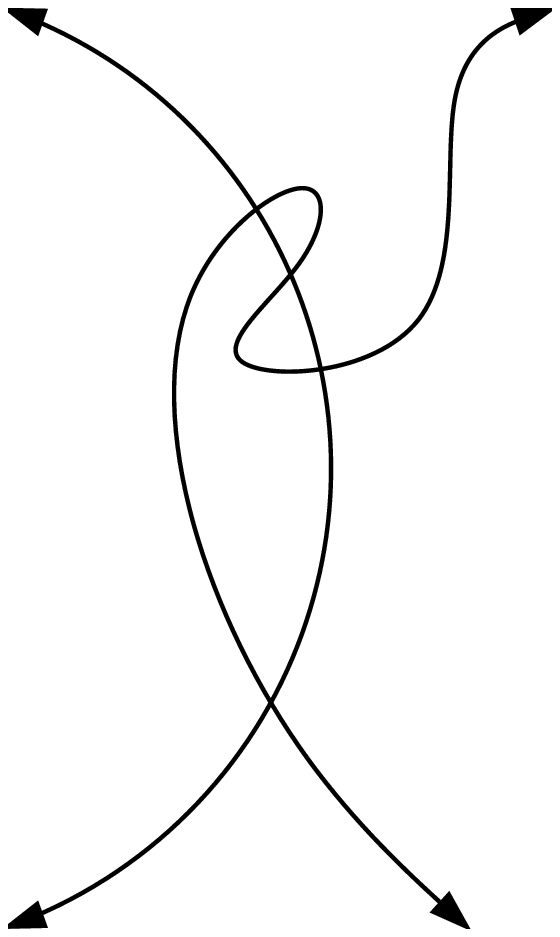}}
    $\Longrightarrow$
    \raisebox{-.75in}{\includegraphics[height=1.5in]{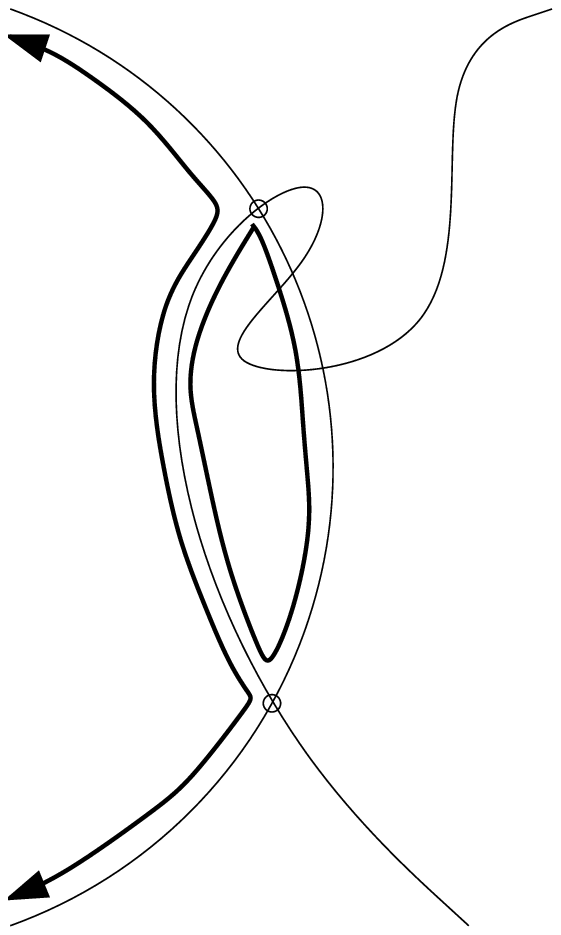}}
    \caption{Compression of a disk}
    \label{fig:compression}
  \end{figure}
  the left image shows the intersection of $B_1$ with $B_2$, while the
  right image shows the constructed $B'_2$ and $B''_2$ in thick lines.
  Since $B'_2$ intersects $B_1$ less than does $B_2$, then it is
  descent equivalent to it.  And since $B'_2$ does not intersect $B_2$
  at all, it is descent equivalent to it.  So both $B_1$ and $B_2$ are
  descent equivalent to $B'_2$, which is a contradiction.
  
  We assume from now on that all components of the intersection of
  $B_1, B_2$ are of the remaining types.  Then $B_1, B_2$ are vertical
  bands and not spheres or disks.  Furthermore their boundaries are
  essential circles in $\boundary M$, otherwise, we could cap them to
  form disks or spheres.
  
  \emph{Case a3 and m3}: The curve is a circle in $B_1$ such that it
  and one component of $\boundary B_1$ bound an annulus $A$.  Then we
  can find a curve $C$ such that $A$ does not contain any other
  component of the intersection.  Then $C$ divides $B_2$ into an
  annulus $A'$ and another component $B'$, which is either another
  annulus or a M\"obius strip.  Indeed, $B_2$ must be of type either
  a3 or m3 since if it was of type m5, then $B_1$ would also be of
  type m5.  One of $A'$ or $B'$ in union with $A$ forms a vertical
  band $B$, which is essential since $\boundary B$ is made of circles
  essential in $\boundary M$.  And after isotopy, $B$ has a simpler
  intersection with both $B_1$ and $B_2$ than they do each other, so
  it is descent equivalent to both, which is a contradiction.
  
  \emph{Case m5}: Let $C$ be isotopic to the core of $B_1$, so it is
  also isotopic to the core of $B_2$, which must also be of type m5.
  Then it cuts each into an annulus, and after an isotopy, the two
  together form an essential vertical annulus $A$ with simpler
  intersection with both $B_1$ and $B_2$ than they do each other, so
  it is descent equivalent to both, which is a contradiction.
  
  \emph{Case a4 and m4}: Finally, the intersection consists only of
  vertical arcs, so a regular neighborhood of the union of the two
  essential vertical bands has at its boundary vertical bands, each
  disjoint from either $B_1$ or $B_2$.  If any band is essential, then
  it is descent equivalent to both $B_1$ and $B_2$, which is a
  contradiction, but if none is essential, then one of them separates
  $B_1$ and $B_2$ from the link, which contradicts the fact that they
  are essential.
\end{proof}

\section{Twisted Link Diagrams}

Theorem~\ref{thm:link_diagrams} will follow from the following lemma.

\begin{lem}
  There is a bijection between each of the following:
  \begin{enumerate}
  \item\label{itm:lot} Ambient isotopy equivalence classes of links in
    stable oriented thickenings.
  \item\label{itm:lsd} Reidemeister equivalence classes of link
    diagrams in stable surfaces.
  \item\label{itm:ald} Reidemeister equivalence classes of abstract
    links.
  \item\label{itm:ld} Reidemeister equivalence classes of twisted link
    diagrams.
  \end{enumerate}
\end{lem}

\begin{figure}[htbp]
  \centering
  \raisebox{-.575in}{\includegraphics[height=1.25in]{Onefoil}}~
  $\Leftrightarrow$
  \raisebox{-.575in}{\includegraphics[height=1.25in]{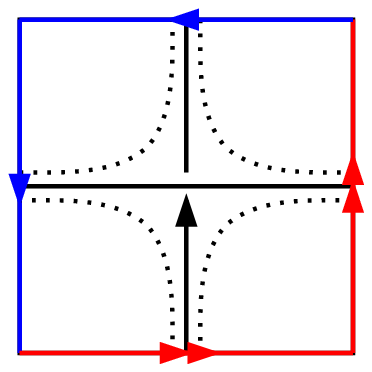}}~
  $\Leftrightarrow$
  \raisebox{-.575in}{\includegraphics[height=1.25in]{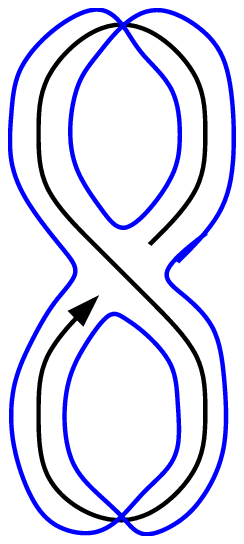}}~
  $\Leftrightarrow$
  \raisebox{-.575in}{\includegraphics[height=1.25in]{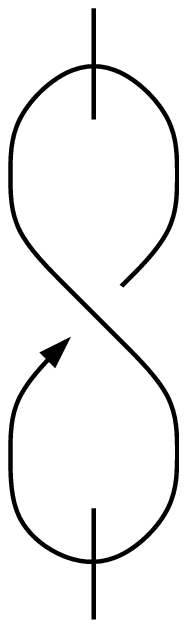}}
  \caption{Creating a diagram for a twisted link}
  \label{fig:linktodiagram}
\end{figure}

Figure~\ref{fig:linktodiagram}
illustrates the process of creating a link diagram from a link in a
thickened Klein bottle.  The first step is a regular projection of the
link to create a link diagram in the zero-section.  The second step
creates an abstract link diagram in a regular neighborhood of the
link diagram in the surface.  The third step immerses this diagram in
the plane, and decorates the result with virtual crossings and bars.

\begin{proof}
  $(\ref{itm:lot})\Leftrightarrow(\ref{itm:lsd})$: Fix an oriented
  thickening of a closed surface and consider the surface as the
  zero-section of the oriented thickening.  By Sard's theorem, the set
  of embeddings of a link in the manifold that have regular
  projections to the surface is dense in the set of all
  ambient-isotopic embeddings that form a link.  The surface
  projection defines a resolution of the double points of the
  projection as follows.  Given the fiber over a neighborhood about
  the double point, if both domain strands are on the same side of the
  neighborhood or on the neighborhood, the image of the strand whose
  doubled point is closer to the boundary of its side has its $\eta$
  vector towards that boundary while the other strand has its
  $\other\eta$ vector towards the other boundary.  And if the strands
  are on different sides, their $\eta$'s are towards their side's
  boundary.  Define a separation of the surface projection from the
  resolution as $\omega = \other\tau\times\eta, \other\omega =
  \tau\times\other\eta$.
    
  On the other hand, given a link surface diagram, choose an oriented
  thickening of the surface.  Define a resolution of the surface
  projection from the separation as $\eta = \omega\times\other\tau,
  \other\eta = \other\omega\times\tau$, and resolve the double point
  by the vectors.  If the thickening with the opposite orientation had
  been chosen, the two links in oriented thickenings are related by an
  orientation-preserving homeomorphism.
    
  By Hudson and Zeeman~\cite{MR29:621}, two embeddings of links in an
  oriented thickening are ambient isotopic if and only if they are
  ambient isotopic by linear moves in arbitrary small neighborhoods.
  Then we may choose neighborhoods of curves such that their
  projections to the surface are generic curves whose intersections
  are in disks.  Then we are in the classical case, and have that
  links that have regular projections are ambient isotopic in the
  manifold if and only if their link surface diagrams are Reidemeister
  equivalent.
    
  There exists a destabilization of the link in oriented thickening
  along a topologically vertical band if and only if there is a
  destabilization of the link surface diagram possibly preceded by a
  sequence of Reidemeister moves.  Indeed, a vertical band exists as
  the fiber over a simple closed path in the projection surface.  And
  if the band is only topologically vertical, it can be made vertical
  through an isotopy of the link that corresponds to a sequence of
  Reidemeister moves because it is isotopic to a vertical band.
  
  $(\ref{itm:lsd})\Leftrightarrow(\ref{itm:ald})$: Given a link
  surface diagram, a regular neighborhood of the link is an abstract
  link diagram.  Given an abstract link diagram, fill in disks along
  its boundary components to obtain a cellular link surface diagram,
  then stabilize the surface to obtain the original link surface
  diagram.
    
  Reidemeister moves on link surface diagrams are done in disk
  neighborhoods obtained after stable isotopy of the link brings the
  strands together.  And Reidemeister moves on abstract link diagrams
  are done in disk neighborhoods obtained by bringing the strands
  together with homeomorphisms of the boundary.
    
  $(\ref{itm:ald})\Leftrightarrow(\ref{itm:ld})$: Obtain an abstract
  link diagram from a twisted link diagram, as follows.  Choose disk
  neighborhoods of the crossings in the plane.  Consider the sphere
  $S^2$ to be $S^2\times0\subset S^2\times I$ and thicken the arcs
  into handles that taper to the boundaries of the disk neighborhoods.
  Figure~\ref{fig:archandles}
  \begin{figure}[htbp]
    \centering
    \includegraphics[height=1.5in]{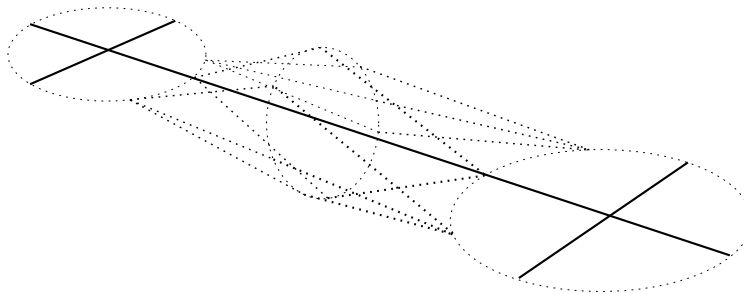}
    \caption{Arc thickened into a tapered handle.}
    \label{fig:archandles}
  \end{figure}
  shows such a tapered handle between two crossings.  Define arc
  neighborhoods in the tapered handles by sliding an interval $I$ from
  one tapered end to the other along arc such that the interval is
  always normal to the arc, and passes through the normal to $S^2$ in
  $S^2\times I$ if and only if there is a bar at that point on the
  arc.  Move the arc neighborhoods by a small amount so that they do
  not intersect one another into one of the two ways shown in
  Figure~\ref{fig:vresolve}.
  \begin{figure}[htbp]
    \centering
    \raisebox{-.45in}{\includegraphics[height=1in]{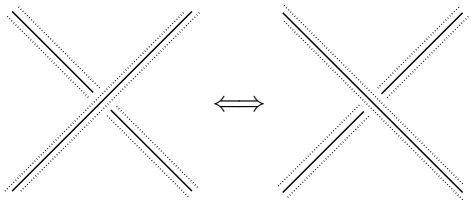}}
    $\Longleftrightarrow$
    \raisebox{-.45in}{\includegraphics[height=1in,angle=90,origin=c]{nbd1v}}
    \caption{Two ways to resolve a virtual crossing}
    \label{fig:vresolve}
  \end{figure}
  Different choices of movement of the arc neighborhoods and of
  direction of turn of the interval as it slides along the arc yield
  homeomorphic abstract link diagrams although it produces a different
  embedding of that diagram in $S^2\times I$.
    
  Suppose two link diagrams differ by an extended Reidemeister move.
  If the move is classical, then the equivalent classical move
  exists on the abstract link diagram.  If the move is virtual or
  twisted, then the embedding is changed in some version of what is
  shown in Figure~\ref{fig:virtualchange}
  \begin{figure}[htbp]
    \centering
    \raisebox{-.35in}{\includegraphics[height=.8in]{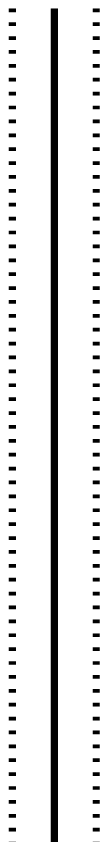}}   
    $\buildrel\text{V1}\over\Longleftrightarrow$
    \raisebox{-.35in}{\includegraphics[height=.8in]{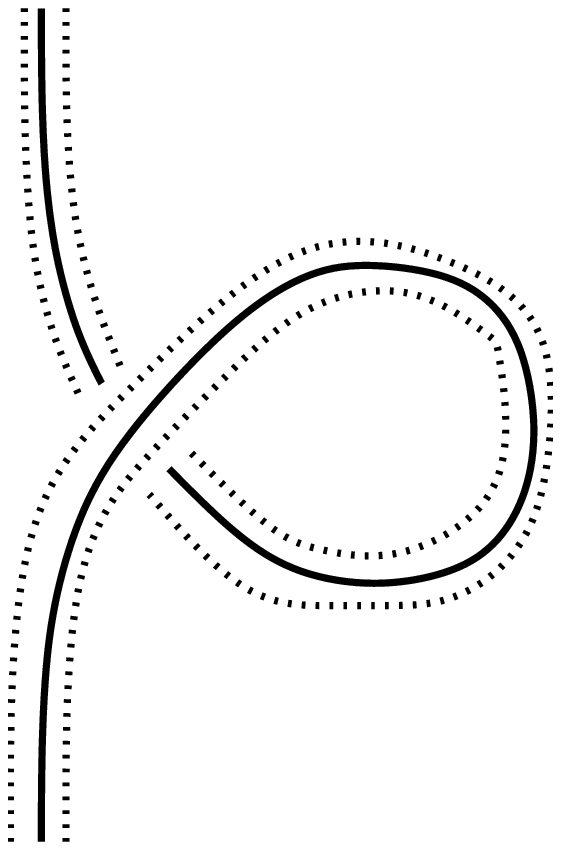}}
    \hspace{1ex}
    \raisebox{-.35in}{\includegraphics[height=.8in]{nbd2vert}}
    $\buildrel\text{V2}\over\Longleftrightarrow$
    \raisebox{-.35in}{\includegraphics[height=.8in]{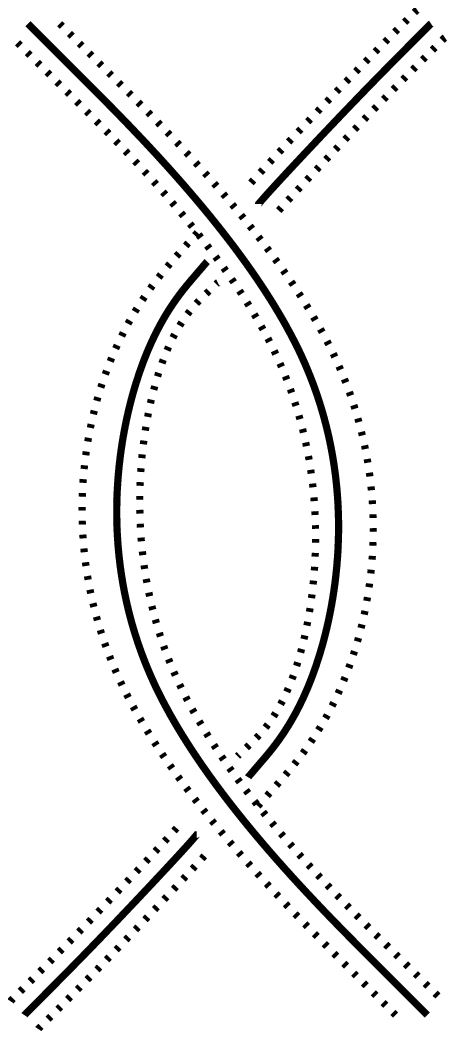}}
    \hspace{1ex}
    \raisebox{-.35in}{\includegraphics[height=.8in]{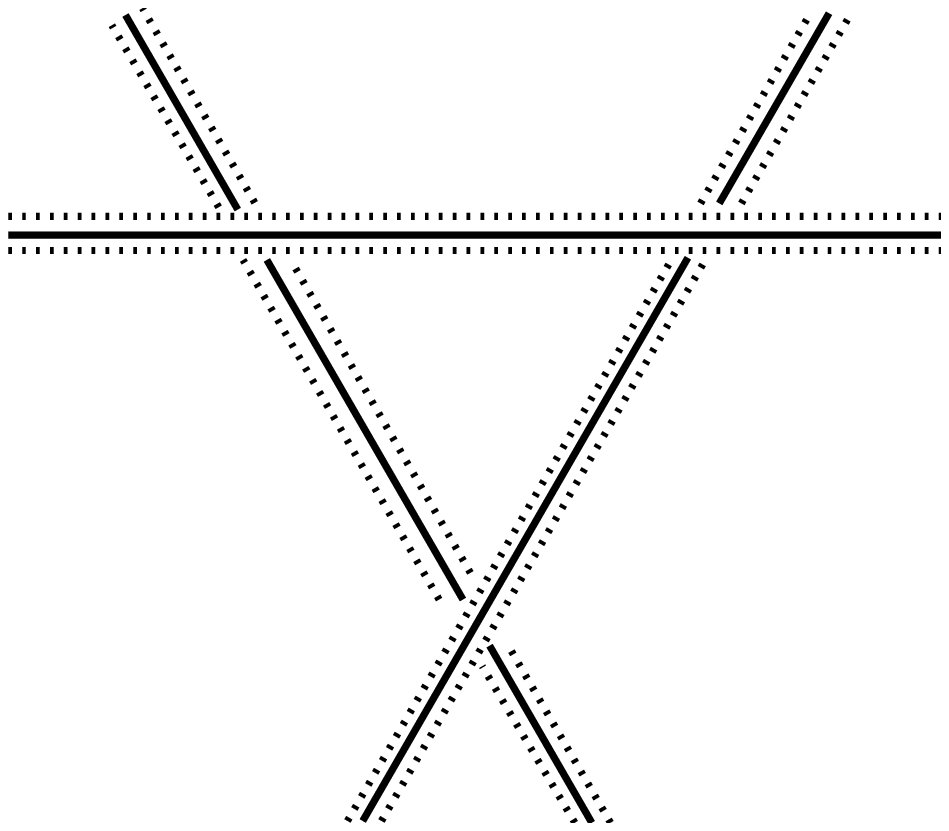}}
    $\buildrel\text{V3}\over\Longleftrightarrow$
    \raisebox{-.35in}{\includegraphics[height=.8in,angle=180,origin=c]{nbd3v}}
    \hspace{1ex}
    \raisebox{-.35in}{\includegraphics[height=.8in]{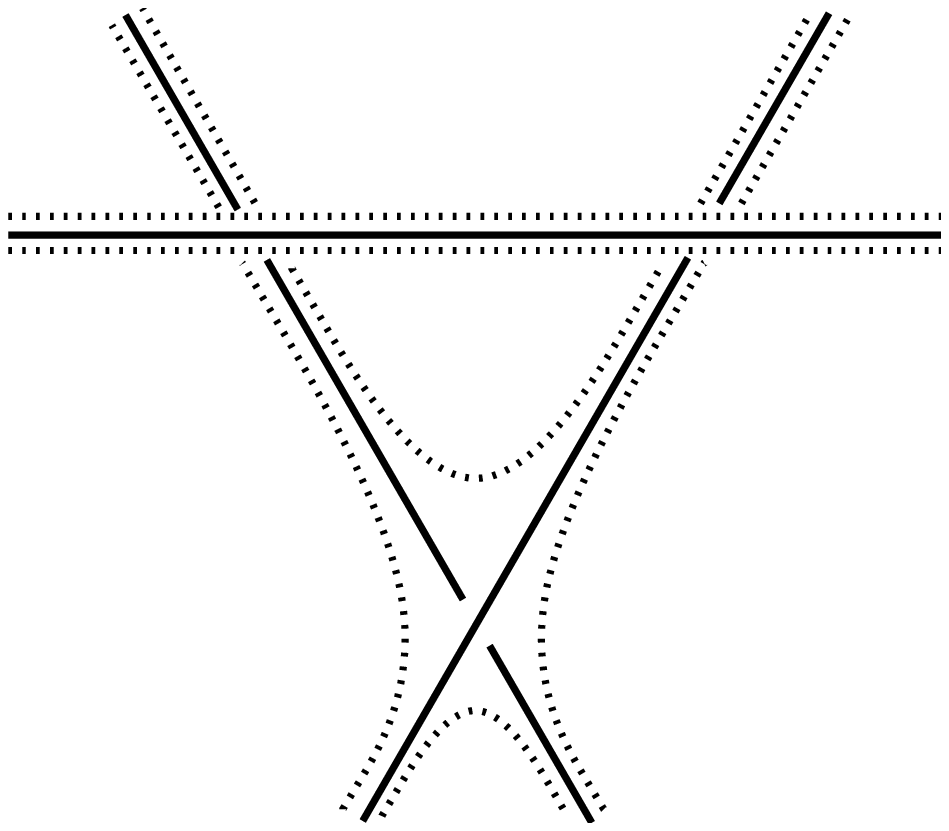}}
    $\buildrel\text{V4}\over\Longleftrightarrow$
    \raisebox{-.35in}{\includegraphics[height=.8in,angle=180,origin=c]{nbd1lr2v}} \\[1ex]
    \raisebox{-.35in}{\includegraphics[height=.8in]{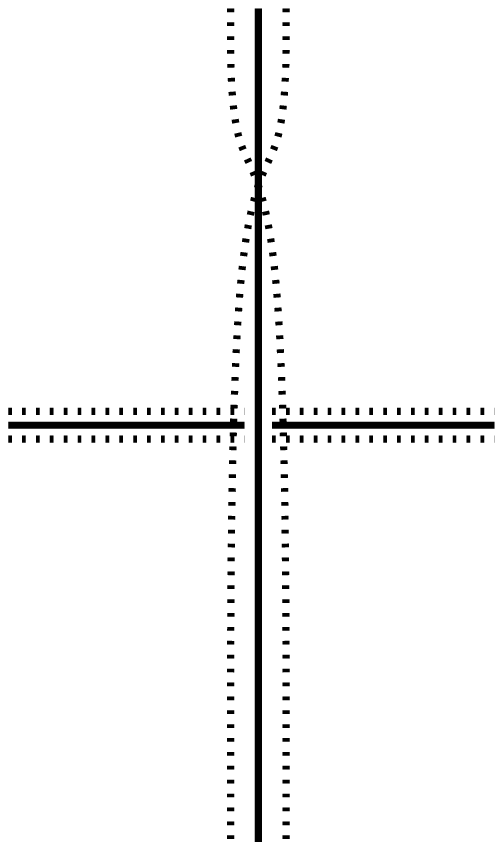}}   
    $\buildrel\text{T1}\over\Longleftrightarrow$
    \raisebox{-.35in}{\includegraphics[height=.8in,angle=180,origin=c]{nbd1vert1v1b}}
    \hspace{1em}
    \raisebox{-.35in}{\includegraphics[height=.8in]{nbd1vert}}
    $\buildrel\text{T2}\over\Longleftrightarrow$
    \raisebox{-.35in}{\includegraphics[height=.8in]{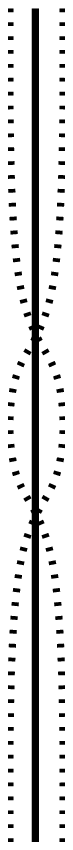}}
    \hspace{1em}
    \raisebox{-.35in}{\includegraphics[height=.8in]{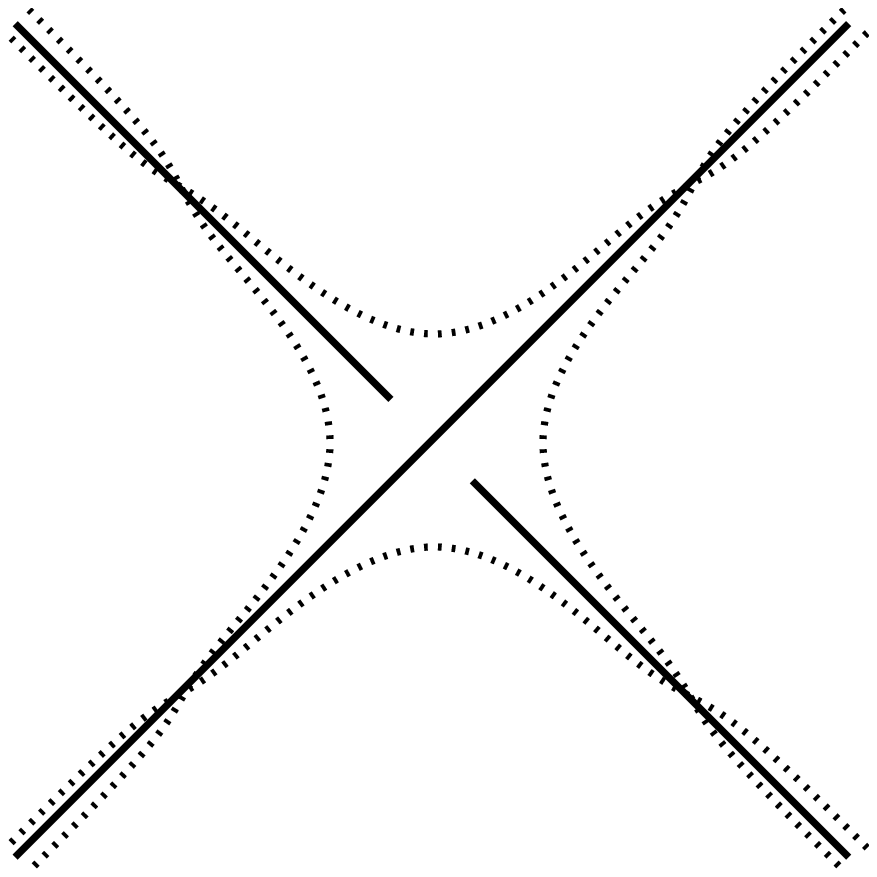}}
    $\buildrel\text{T3}\over\Longleftrightarrow$
    \raisebox{-.35in}{\includegraphics[height=.8in]{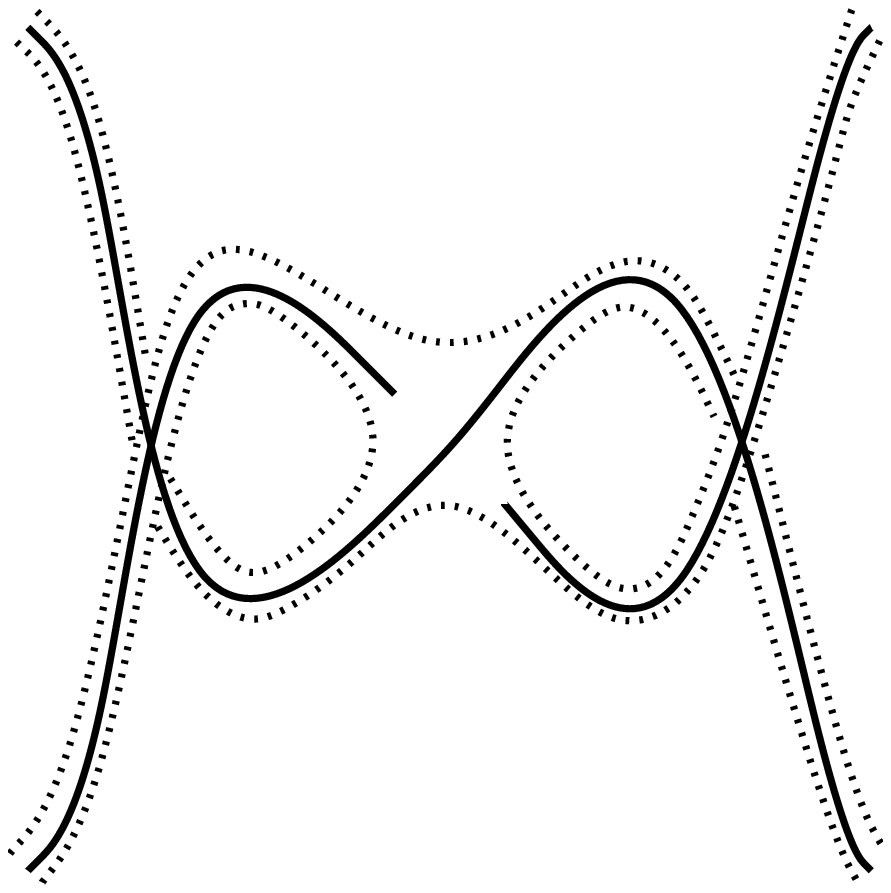}}
    \caption{Virtual and twisted changes to embeddings.}
    \label{fig:virtualchange}
  \end{figure}
  but the abstract link diagram is unchanged.
    
  Obtain a twisted link diagram from an abstract link diagram as
  follows.  Embed the abstract link diagram in $S^2\times I$ such that
  each crossing has a disk neighborhood in $S^2$, and that the
  projection of the embedding to $S^2$ is an immersion.  On the
  projected diagram, draw crossings in the over-under form depending
  on the resolution of the abstract link crossings.  Draw each
  intersection of two projected arcs as a virtual crossing.  Draw each
  inter-component intersection of the projection of the two boundary
  components of the neighborhood of each arc as bars on the arc.
  Figure~\ref{fig:arcbend}
  \begin{figure}[htbp]
    \centering
    \includegraphics[height=1.5in]{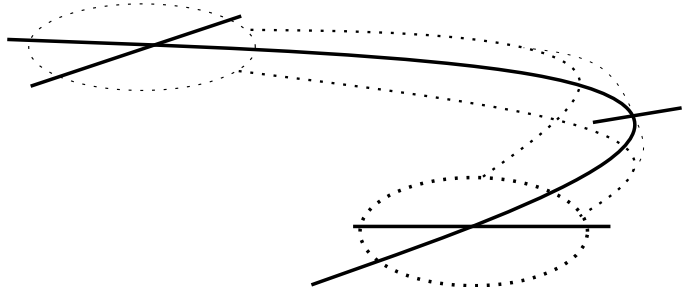}
    \caption{Arc boundary intersection giving rise to a bar}
    \label{fig:arcbend}
  \end{figure}
  shows the projection of an arc and the boundary components, and
  draws one bar on the arc for the single inter-component intersection
  of thee two boundary components.
    
  Suppose two different embeddings $g,g'$ of the same abstract link
  give rise to different twisted link diagrams.  We can change $g$ by
  embedding each neighborhood of a crossing with the opposite rotation
  of the arcs to make them match that of $g'$, which corresponds to a
  series of T3 moves.  Figure~\ref{fig:flipping}
  \begin{figure}[htbp]
    \centering
    \raisebox{-.5in}{\includegraphics[height=1in]{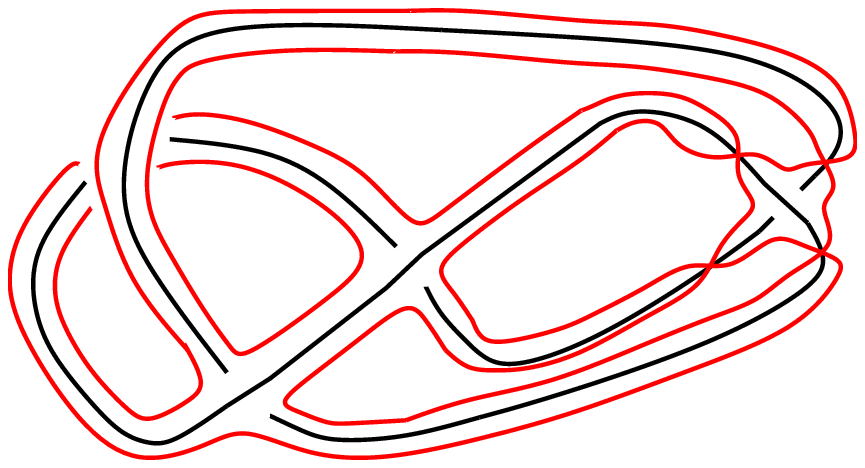}}~
    $\Leftrightarrow$~
    \raisebox{-.5in}{\includegraphics[height=1in]{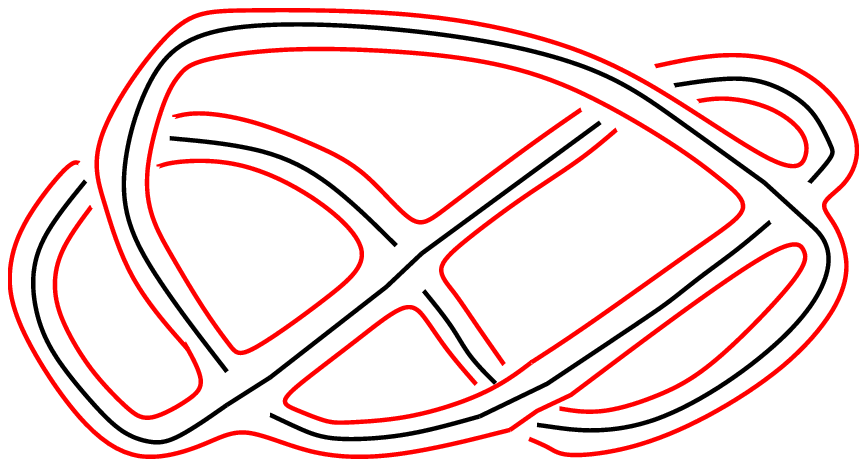}}
    \raisebox{-.5in}{\includegraphics[height=1in]{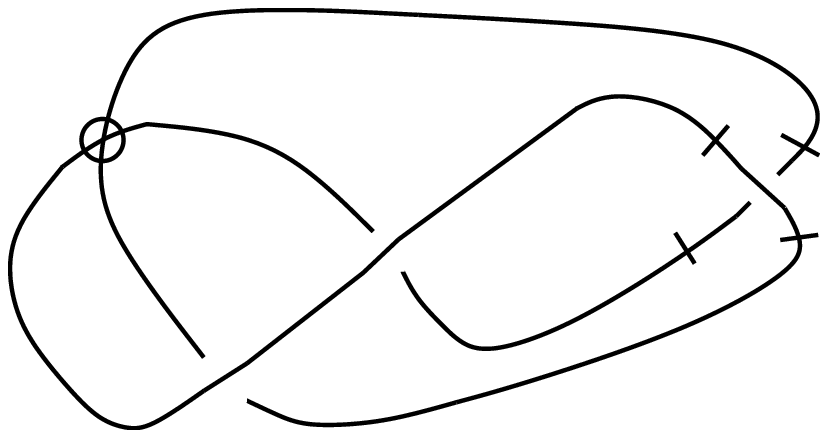}}~
    $\Leftrightarrow$~
    \raisebox{-.5in}{\includegraphics[height=1in]{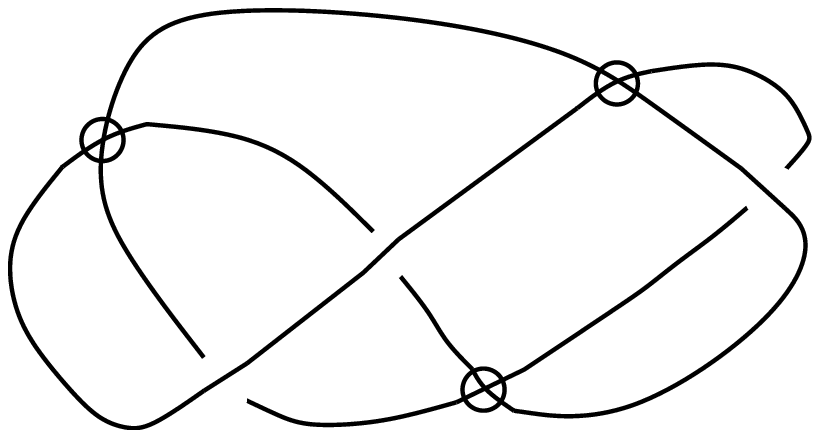}}
    \caption{Performing a T3 move}
    \label{fig:flipping}
  \end{figure}
  shows the effect of performing a T3 move on both an embedding of an
  abstract link and on the corresponding twisted link diagram.  We can
  change $g$ by an isotopy of $S^2$ to make its neighborhoods of
  crossings including the arcs correspond to that of $g'$, which
  corresponds to a series of V4 moves.  For each arc neighborhood in
  turn, we can change $g$ by a homotopy until the arc coincides with
  that of $g'$, which corresponds to a finite sequence of V1--4 and T1
  moves, and we can then change $g$'s embedding of the arc
  neighborhood itself so that it matches that of $g'$, which
  corresponds to a finite sequence of T1--2 moves.
    
  Suppose two abstract links differ by a single R1--3 move.  By the
  above, we can choose embeddings of these links that coincide outside
  of a disk neighborhood containing the crossings and edges involved.
  Then the twisted link diagrams will differ by the corresponding
  R1--3 move.
\end{proof}

Since the projection of a link is a generic immersion of curves in the
surface homeomorphic to the zero-section, then projections of a link
in an oriented thickening are classified by intersigned Gauss codes.
By adding a writhe sign to an intersigned Gauss code and considering
it up to abstract Reidemeister equivalence defined in
Figure~\ref{fig:reidilc}, we then obtain the following:

\begin{cor}
  Reidemeister equivalence classes of link diagrams in stable surfaces
  correspond to Reidemeister equivalence classes of intersigned link
  codes.
\end{cor}

\section{The Twisted Jones Polynomial is an Invariant}

\begin{figure}[htbp]
  \centering
  \raisebox{-.5in}{\includegraphics[height=1in]{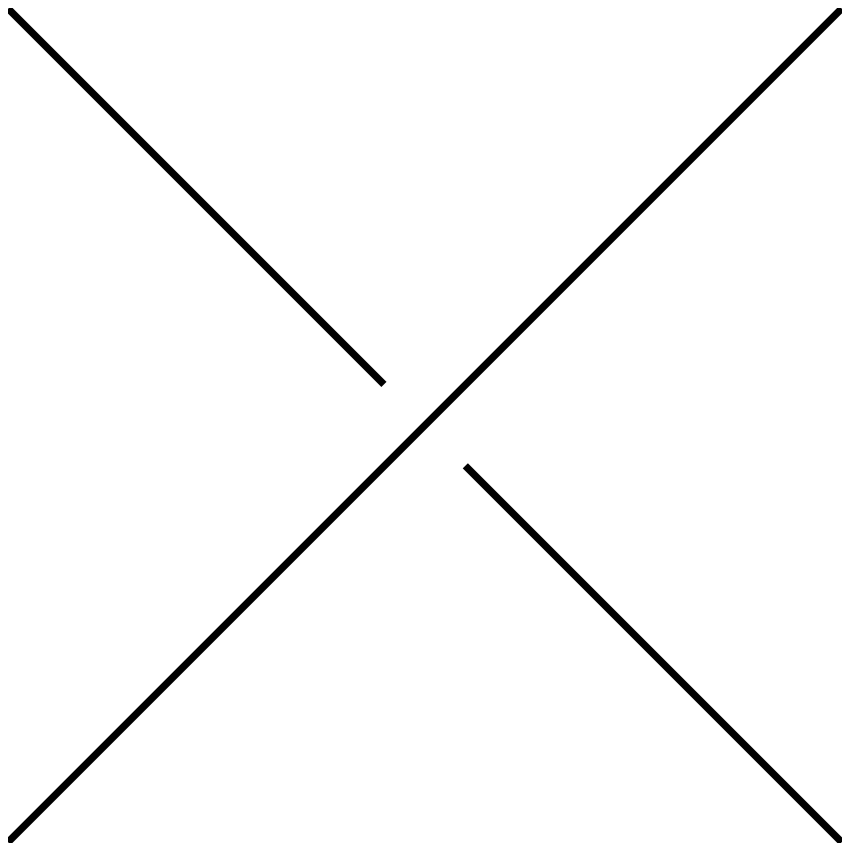}}~
  $\buildrel a\over\Rightarrow$~
  \raisebox{-.5in}{\includegraphics[height=1in]{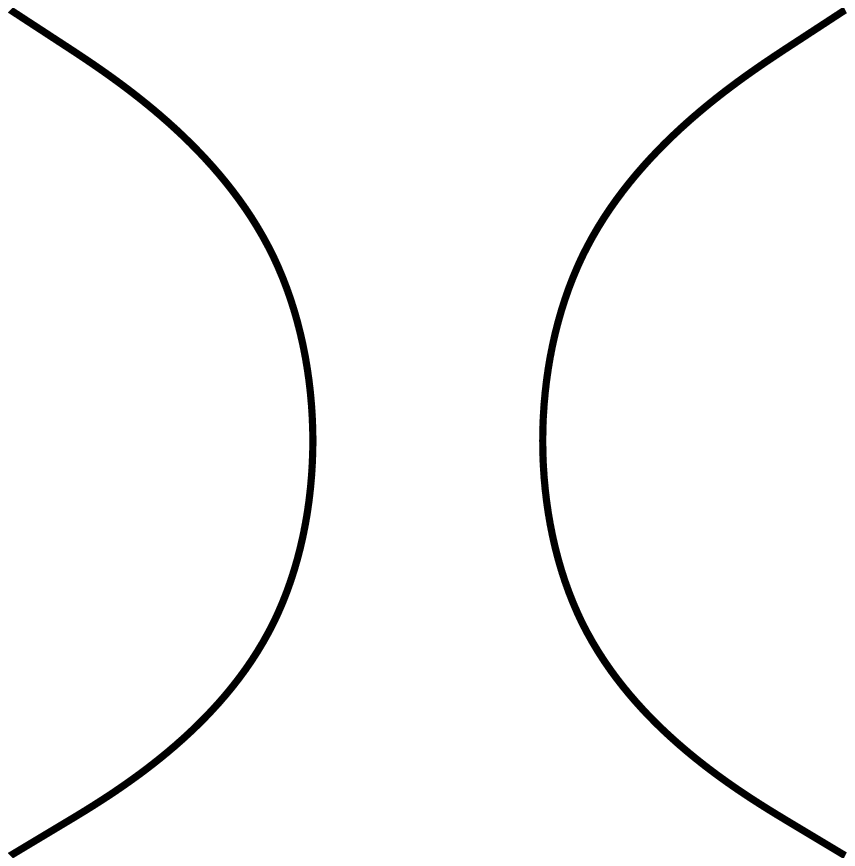}}~
  or~
  \raisebox{-.5in}{\includegraphics[height=1in]{lr}}~
  $\buildrel b\over\Rightarrow$~
  \raisebox{-.5in}{\includegraphics[height=1in,angle=90,origin=c]{lrA}}
  \caption{Two possible smoothings of a crossing.}
  \label{fig:smoothings}
\end{figure}

A \emph{smoothing} of a real crossing of a link diagram is a
transformation of one of two kinds labeled ``a'' and ``b'' as shown in
Figure~\ref{fig:smoothings}.
A \emph{state} of a link diagram is a collection of smoothings of all
its crossings.  The \emph{bracket polynomial} of a link diagram $D$ is
an element of $\rZ[A^{\pm1},M]$ defined from the possible states
$\mathcal S(D)$ of $D$ by:
\begin{align*}
  \left<D\right> &= \sum_{S\in\mathcal
    S(D)}A^{a(S)-b(S)}(-A^{-2}-A^2)^{c(S)}M^{d(S)}
\end{align*}
where:
\begin{itemize}
\item $a(S)$ is the number of $a$-smoothings,
\item $b(S)$ is the number of $b$-smoothings,
\item $c(S)$ is the number of circles with an even number of bars, and
\item $d(S)$ is the number of circles with an odd number of bars.
\end{itemize}
The \emph{twisted Jones polynomial} of a link diagram $D$ is an
element of $\rZ[A^{\pm1},M]$ calculated as:
\begin{align*}
  V_D(A,M) = (-A)^{-3w(D)}\left<D\right>.
\end{align*}

We prove that the extension to the Jones polynomial defined above is
an invariant of links in oriented thickenings, and that it
distinguishes a class of links in oriented thickenings from virtual
links.

\begin{proof}[Proof of Theorem~\ref{thm:Jones_polynomial}]
  The bracket polynomial of link diagrams in oriented thickenings is
  equivalent to that defined by the relations:
  \begin{enumerate}
  \item $\left<\emptyset\right> = 1$,
  \item\label{itm:smooth}
    $\left<\raisebox{-.155in}{\includegraphics[height=.4in]{lr}}\right>
    =
    A\left<\raisebox{-.155in}{\includegraphics[height=.4in]{lrA}}\right>
    + \inv
    A\left<\raisebox{-.155in}{\includegraphics[height=.4in,angle=90,origin=c]{lrA}}\right>$,
  \item
    $\left<D\disjoint\raisebox{-.14in}{\includegraphics[height=.35in]{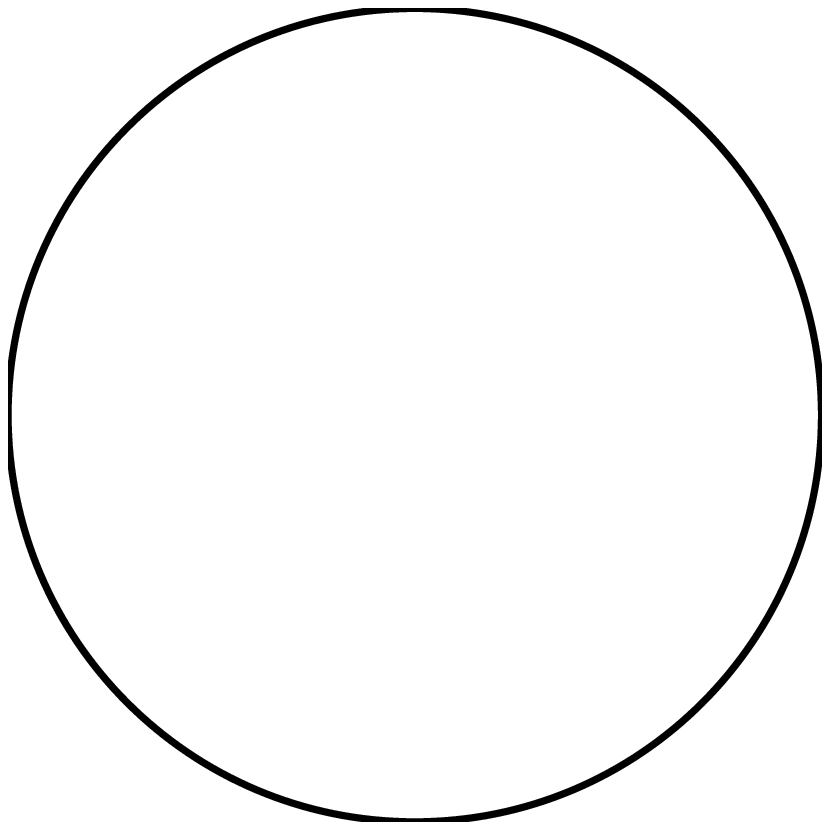}}\right>
    = (-A^{-2}-A^2)\left<D\right>$,
  \item
    $\left<D\disjoint\raisebox{-.155in}{\includegraphics[height=.4in]{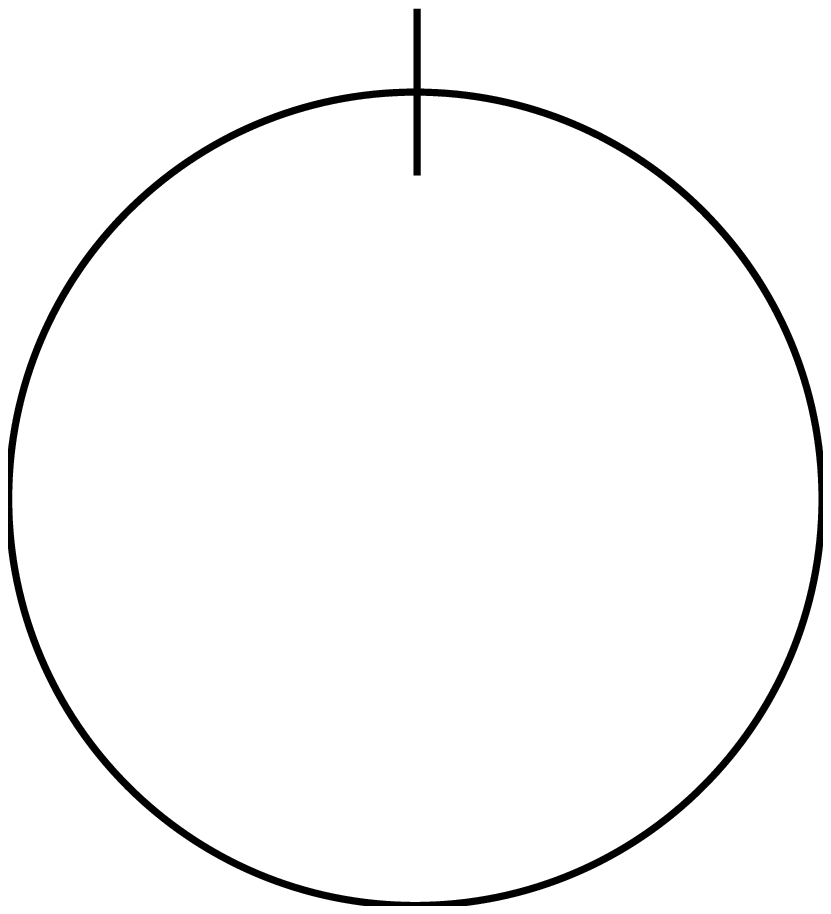}}\right>
    = M\left<D\right>$,
  \item\label{itm:evenbars}
    $\left<\raisebox{-.155in}{~\includegraphics[height=.4in]{1vert}~}\right>
    =
    \left<\raisebox{-.155in}{~\includegraphics[height=.4in]{1vert2bar}~}\right>$.
  \end{enumerate}
  This can be seen by first using relation~\ref{itm:smooth} to smooth
  the crossings of a diagram $D$ in all possible ways, so that each
  final diagram of circles corresponding to one of the states of $D$.
  Then, relation~\ref{itm:evenbars} allows each circle to be reduced
  to either 0 or 1 bars.  Finally, the other relations allow the
  calculation of the bracket polynomial.

  Invariance with respect to R1--3 and V1--4 is immediately due to
  that of the Jones polynomial.  Invariance with respect to T1 is
  immediate since it involves no real crossings.  Invariance with
  respect to T2 is shown by relation~\ref{itm:evenbars}.  Then, the
  twisted Jones polynomial is invariant with respect to move T3 by the
  calculation:
  \begin{align*}
    \left<\raisebox{-.175in}{\includegraphics[height=.4in]{1lr4bar}}\right>
    &=
    A\left<\raisebox{-.175in}{\includegraphics[height=.4in]{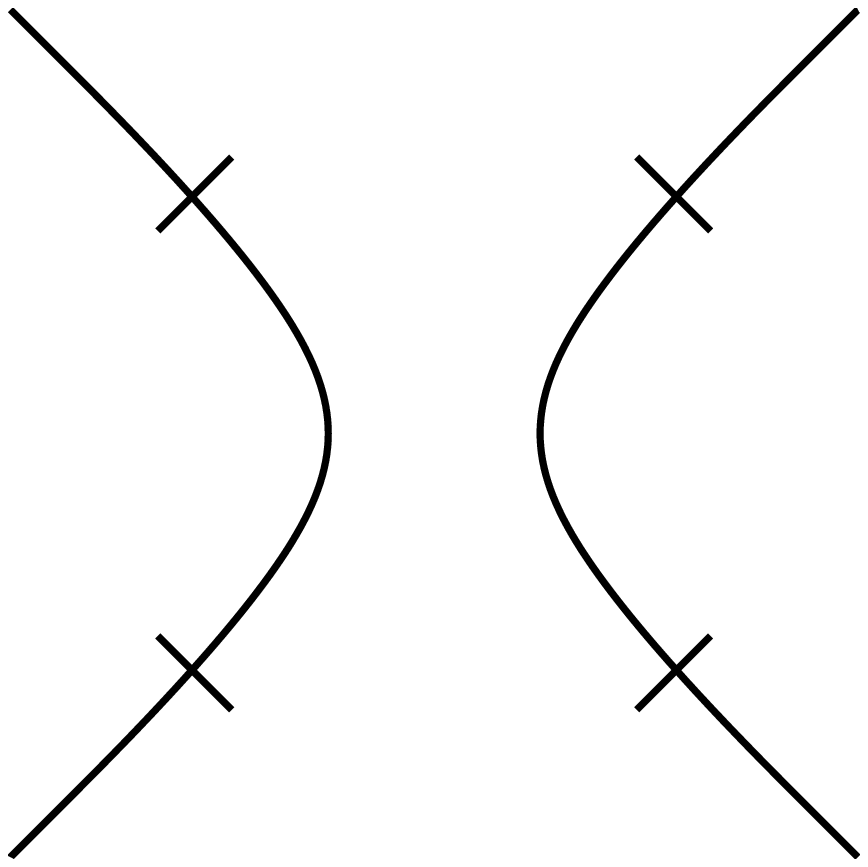}}\right>+
    \inv
    A\left<\raisebox{-.175in}{\includegraphics[height=.4in,angle=90,origin=c]{1lr4barA}}\right>
    &\text{by (\ref{itm:smooth})} \\
    &=
    A\left<\raisebox{-.175in}{\includegraphics[height=.4in]{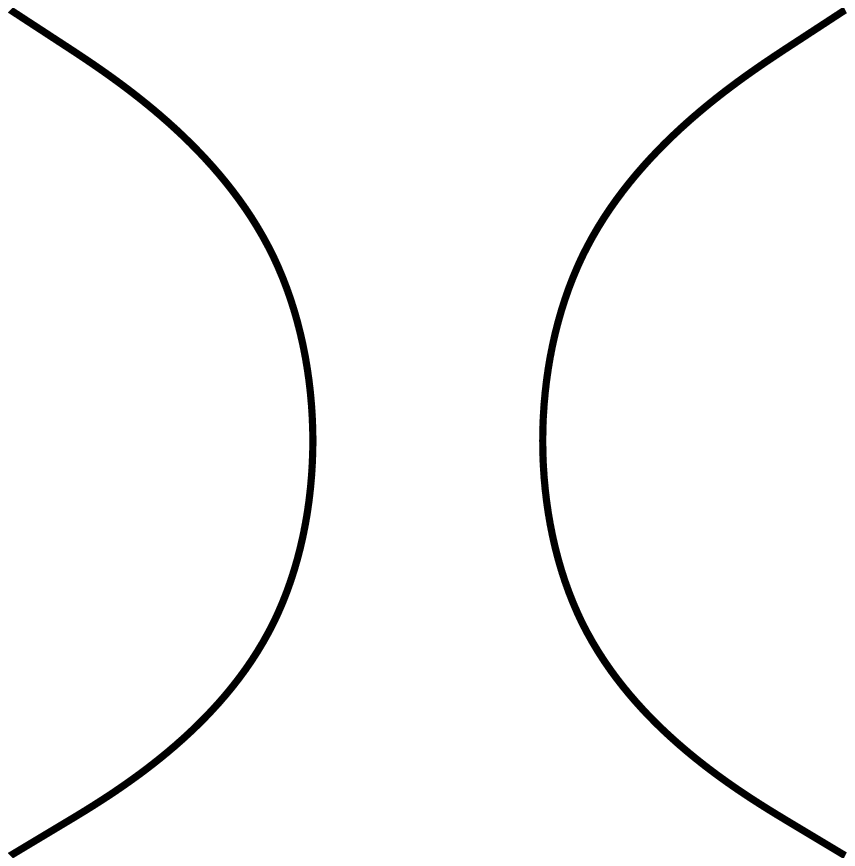}}\right>+
    \inv
    A\left<\raisebox{-.175in}{\includegraphics[height=.4in,angle=90,origin=c]{1lrA}}\right>
    &\text{by (T1)} \\
    &=
    A\left<\raisebox{-.175in}{\includegraphics[height=.4in]{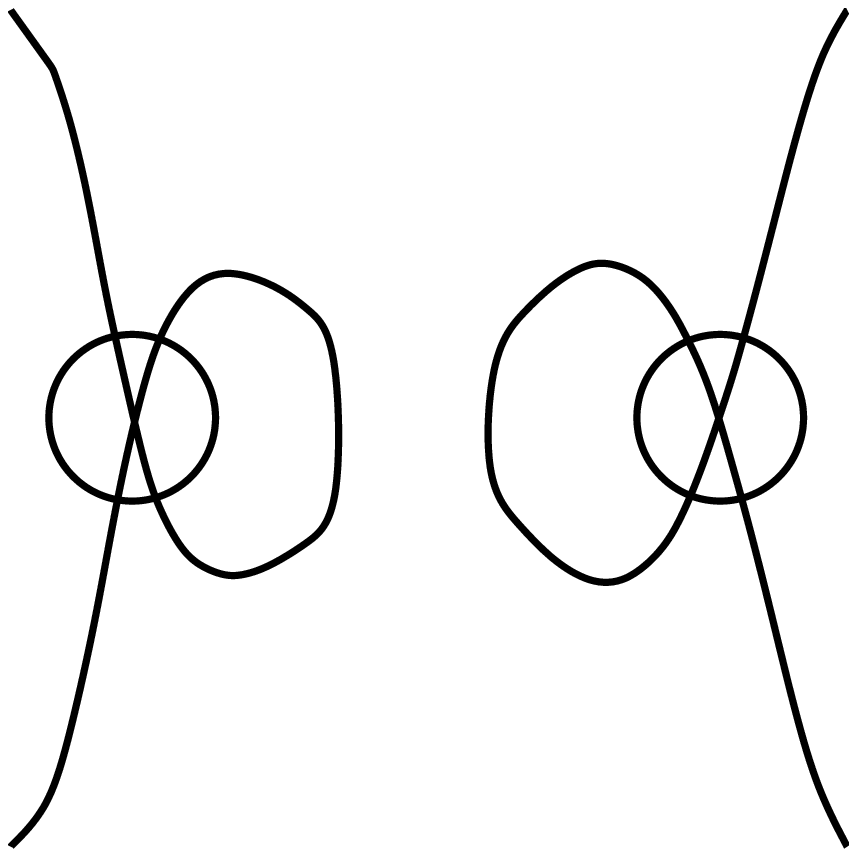}}\right>+
    \inv
    A\left<\raisebox{-.175in}{\includegraphics[height=.4in]{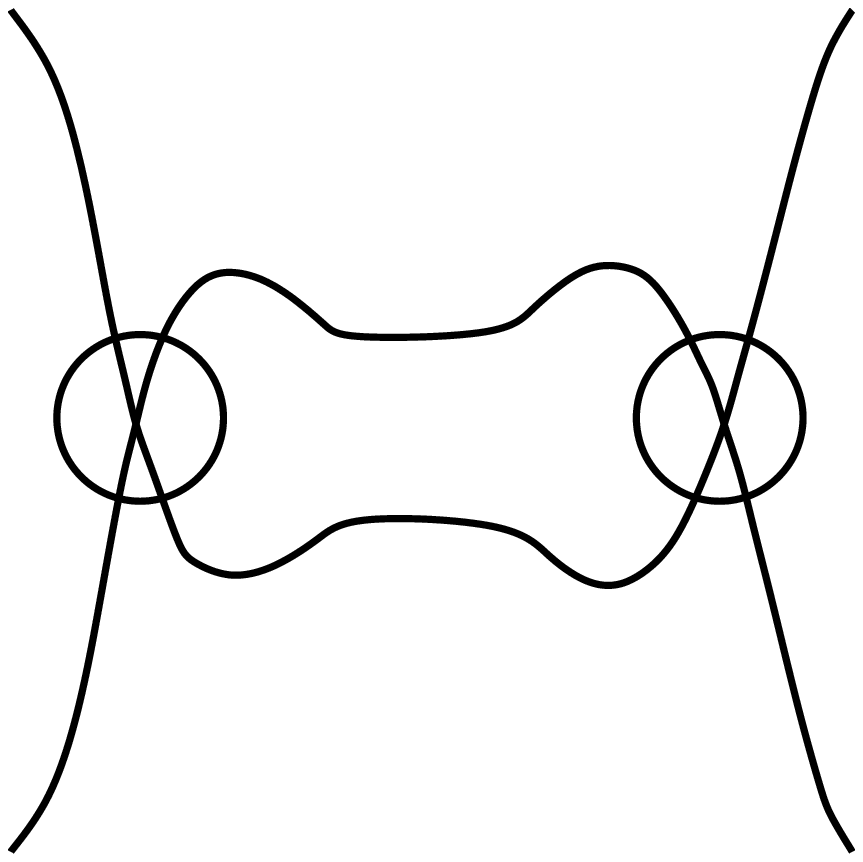}}\right>
    &\text{by (V1,V2)} \\
    &=
    \left<\raisebox{-.175in}{\includegraphics[height=.4in]{1lr2virtflank}}\right>.
    &\text{by (\ref{itm:smooth})}
  \end{align*}
  So the twisted Jones polynomial is an invariant of links in oriented
  thickenings.
  
  Since virtual links do not have bars on their edges, then their
  states will not have any circles with an odd number of bars, and so
  their twisted Jones polynomials will be in $\rZ[A^{\pm1}]$.  And
  since the $(-A^{-2}-A^2)$ term is raised to the number of circles in
  a state, the twisted Jones polynomial will factor into a product of
  $(-A^{-2}-A^2)$ and the Jones polynomial of the diagram.
\end{proof}

The Jones polynomial can be extended to links in oriented thickenings
by the device of ignoring bars on edges.  This invariant corresponds
to setting $M = -A^{-2}-A^2$ in the new polynomial and then dividing the
result by $-A^{-2}-A^2$.

The onefoil shown in Figure~\ref{fig:onefoil} has twisted Jones
polynomial:
\begin{align*}
  V_\text{Onefoil}(A,M) &= A^{-6}+(1-M^2)A^{-2}
\end{align*}
so it is not a virtual knot.

\section{Twisted Jones of Two-Colorable Diagrams}

The knot $K$ in Figure~\ref{fig:T1212} has twisted Jones polynomial:
\begin{align*}
  V_K(A,M) &= (-A^{-2}-A^2)(A^{-4}+A^{-6}-A^{-10})
\end{align*}
which is a product of $(-A^{-2}-A^2)$ and of the Jones polynomial for
the knot.  But this knot is in the projective plane and cannot be
moved to an affine subset of this space, so it is not a classical
knot~\cite{MR91i:57001}.  We will see that the knot's twisted Jones
polynomial factors in this manner because it has a two-colorable
diagram.

\begin{proof}[Proof of Theorem~\ref{thm:link_two-colorable}]
  Let $D$ be a two-colored diagram for a link.  This two-coloring
  partitions the faces of $D$ into two sets of circles each of one of
  the two colors and such that each crossing of the diagram will
  separate a disk neighborhood of that crossing into two pairs of
  regions that are opposite one another with respect to the crossing
  and that contain like-colored parts of the faces of the diagram.
  Figure~\ref{fig:two-colored}
  \begin{figure}[htbp]
    \centering 
    \includegraphics[height=1.5in]{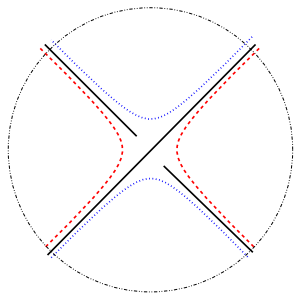}
    \caption{A crossing of a two-colored diagram}
    \label{fig:two-colored}
  \end{figure}
  shows this situation for an arbitrary crossing.  In the figure, the
  parts of faces are drawn either dotted or dashed to represent the
  two colors while the boundary of the disk is drawn with a
  dash-dot-dotted pattern.  The two possible smoothings of each
  crossing of a diagram also pair opposite regions of the crossing.
  In the figure, the dashed lines are the ``$a$'' smoothing while the
  dotted lines are the ``$b$'' smoothing.  Then, the faces of one
  color of $D$ are the circles of the state of $D$ whose smoothings
  pair the corresponding two opposite regions of each crossing, and
  the faces of the other color are the circles of the complementary
  state of $D$.  Since the circles of each of these states correspond
  to faces of a diagram and that by the Jordan curve theorem the faces
  cross the link diagram an even number of times at bars, then the
  circles have an even number of bars.
  
  All other pairs of states of $D$ can be obtained from the above pair
  of states by taking the opposite smoothings at some subset of the
  crossings.  Changing a state at a crossing corresponds to putting
  bars on all the edges coming into that crossing, and the effect of
  this change on the faces is show in
  Figure~\ref{fig:changesmoothing}.
  \begin{figure}[htbp]
    \centering
    \raisebox{-.75in}{\includegraphics[height=1.5in]{lrsmooth}}
    \hspace*{1em}$\Longleftrightarrow$\hspace*{1em}
    \raisebox{-.75in}{\includegraphics[height=1.5in]{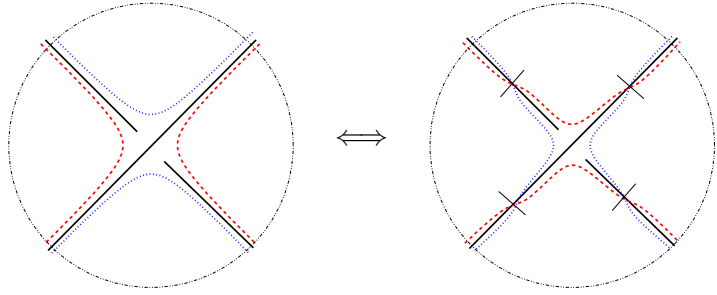}}
    \caption{Changing the smoothings at a crossing.}
    \label{fig:changesmoothing}
  \end{figure}
  In the figure, the diagram on the left shows ``$a$'' smoothing as
  the dashed lines pairing the left and right regions, while the
  ``$b$'' smoothing is the dotted line pairing the top and bottom
  regions, and the diagram on the right has the opposite pairings.
  While adding the bars changes the paths that form the circles of the
  new state, these circles still have an even number of bars.
  Therefore all the states of $D$ have circles with an even number of
  bars on the circles.  Then, the twisted Jones polynomial for this
  link will not have any $M$ term, and will be divisible by
  $-A^{-2}-A^2$.
\end{proof}

\section{The Twisted Link Group is an Invariant}

A \emph{Wirtinger presentation} of a group is a finite presentation
whose relators are of the form $\inv a\inv cbc$ for any three
generators $a,b,c$ of the presentation.  When $c = b$, the relation
becomes the identification $\inv a b$, so identification relations are
admissible in a Wirtinger presentation.

\begin{figure}[htbp]
  \centering
  \includegraphics[height=1.3in]{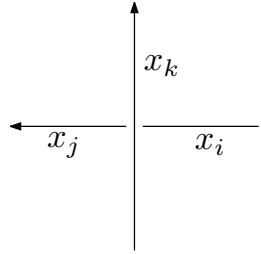}
  \caption{Generators at a crossing for a formal presentation}
  \label{fig:Wcrossing}
\end{figure}

The \emph{group of a link diagram} $\Pi L$ is a group derived from a
link diagram $L$, and has a Wirtinger presentation that has one
generator for every arc of the diagram between undercrossings, and one
relation for every crossing.  Figure~\ref{fig:Wcrossing}
shows three generators of a formal Wirtinger presentation obtained
from a diagram of a virtual link at a crossing.  These generators have
the relation:
\begin{align*}
  x_j &= \inv x_k x_i x_k.
\end{align*}
The group of a link diagram is also called the \emph{upper group} of
the diagram and may be denoted $\Pi^u L$.  The \emph{lower group}
$\Pi^l L$ of a link diagram has one generator for every arc of the
diagram between overcrossings, and again one relation for every
crossing defined analogously to that of the upper group.

\begin{figure}[htbp]
  \centering
  \includegraphics[height=2in]{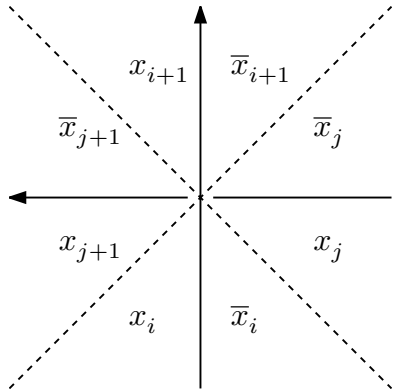}
  \caption{Generators at a crossing for a twisted presentation}
  \label{fig:ncrossing}
\end{figure}

\begin{figure}[htbp]
  \centering
  \includegraphics[height=2in]{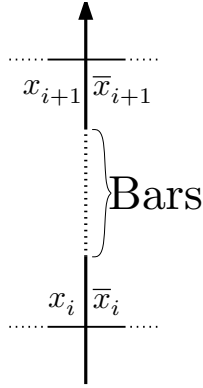}
  \caption{Generators on an edge for a twisted presentation}
  \label{fig:edgegen}
\end{figure}

The \emph{twisted link group of a diagram} $\twisted\Pi L$ is a group
with a Wirtinger presentation that has two generators for each side of
each end of an edge.  Figure~\ref{fig:ncrossing}
shows eight generators obtained from a diagram of a link at a
crossing.  These generators have the four crossing relations:
\begin{align*}
  x_{i+1} &= x_i, & x_{j+1} &= \inv x_i x_j x_i, \\
  {\other x}_{i+1} &= \inv{\other x}_j {\other x}_i {\other
    x}_j, & {\other x}_{j+1} &= {\other x}_j.
\end{align*}
The four generators of an edge have two relations depending on the
parity of the number of bars on the edge.  Figure~\ref{fig:edgegen}
shows the four generators of an edge.  When an edge has an even number
of bars between two crossings, we have the relations:
\begin{align*}
  x_{i+1} &= x_i, & {\other x}_{i+1} &= {\other x}_i, \\
  \intertext{and when there is an odd number of bars on the arc, we
    have the relations:} x_{i+1} &= {\other x}_i, & {\other x}_{i+1} &=
  x_i.
\end{align*}

\begin{proof}[Proof of Theorem~\ref{thm:link_twisted_group}]
  Note that performing moves V1--4 and T1--2 on a link diagram yield a
  link diagram with an identical presentation.  Invariance with
  respect to moves R1--3 is the same as for classical link diagrams
  since neither virtual crossings nor bars are involved.  Finally,
  invariance with respect to T3 is as follows.
  Figure~\ref{fig:T3gens}
  \begin{figure}[htbp]
    \centering
    \raisebox{-.75in}{\includegraphics[height=1.5in]{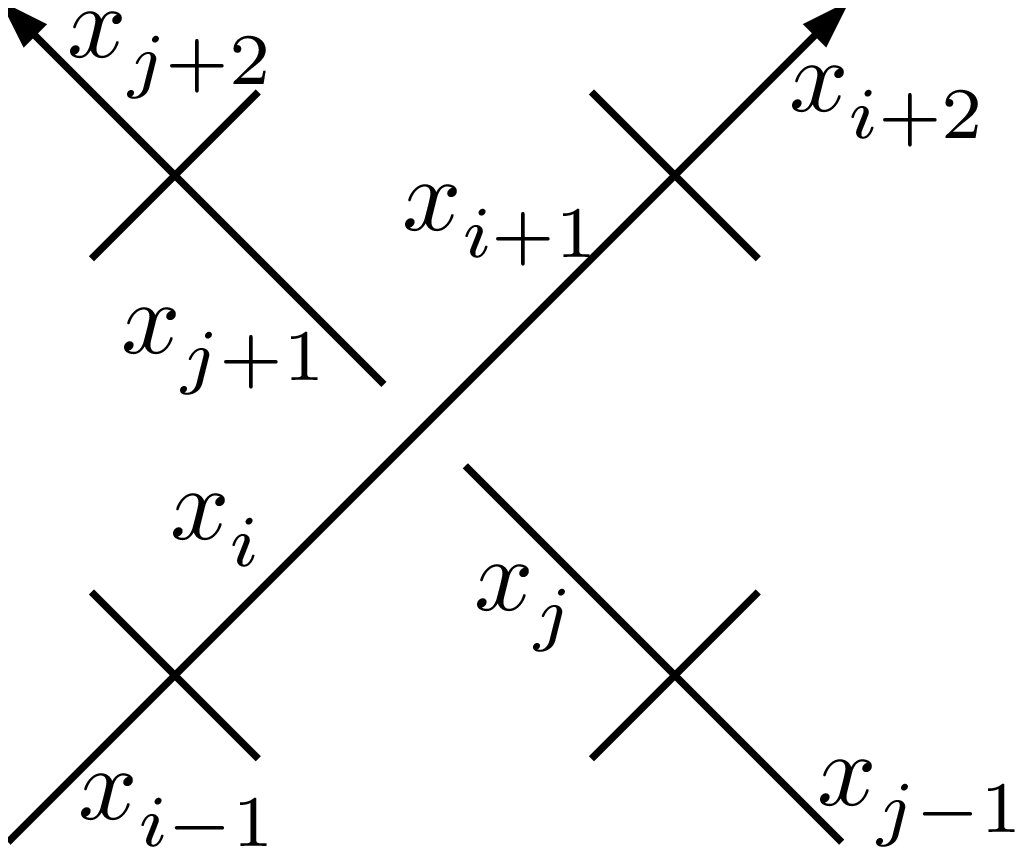}}
    $\buildrel\text{T3}\over\Longleftrightarrow$
    \raisebox{-.75in}{\includegraphics[height=1.5in]{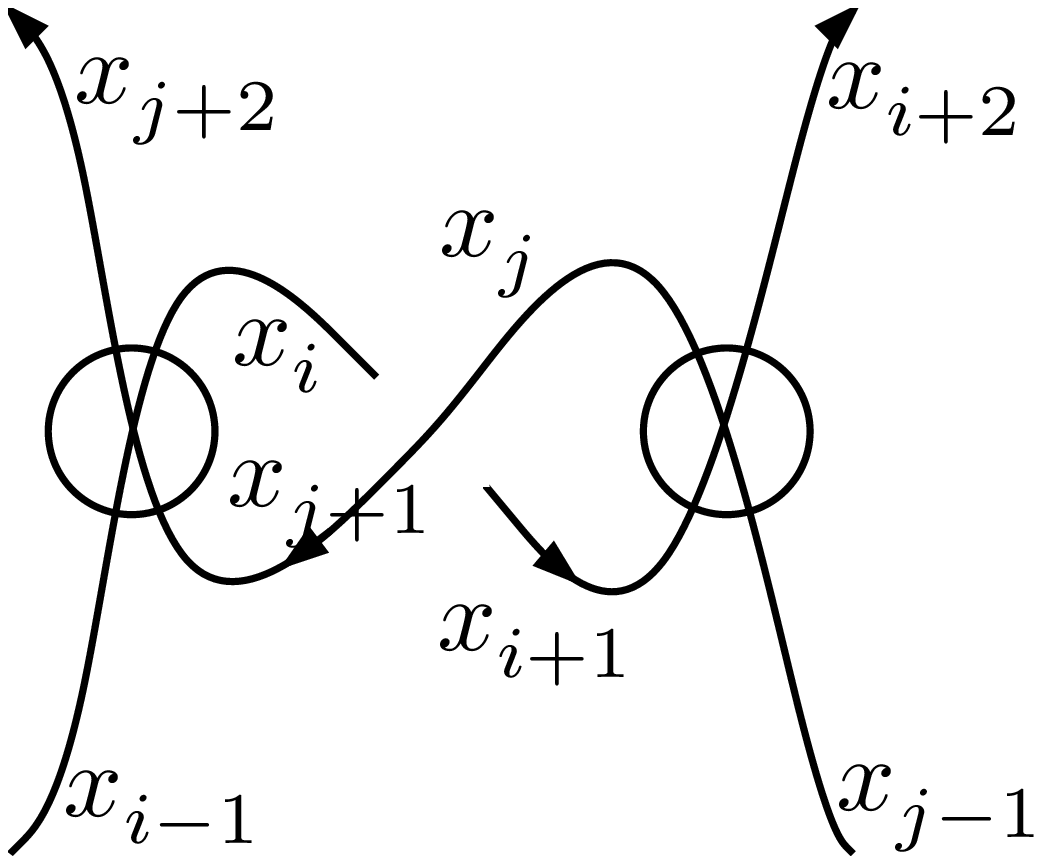}}
    \caption{The generators in move T3.}
    \label{fig:T3gens}
  \end{figure}
  shows half of the generators on the edges of the two diagrams.
  Then, the relations given by the crossing are:
  \begin{align*}
    &\text{On the left} & &\text{On the right} \\
    &\inv x_{j+1} \inv x_i x_j x_i & &\inv x_{j+1} \inv x_i x_j x_i \\
    &\inv x_{i+1} x_i & &\inv x_{i+1} x_i
  \end{align*}
  which are identical.  Doing the same for the other half of the
  generators shows the corresponding equivalence.
\end{proof}


\begin{thebibliography}

\bibitem{bourgoin03:_class_immer_curves}
\textbf{M Bourgoin}, \emph{Classifying Immersed Curves}, in Preparation

\bibitem{bourgoin04:_fundam_group_virtual_knot}
\textbf{M Bourgoin}, \emph{On the Fundamental Group of a Virtual Link}, in
  Preparation

\bibitem{MR91i:57001}
\textbf{Yu\,V Drobotukhina}, \emph{An analogue of the {J}ones polynomial for
  links in {${\bf R}{\rm P}\sp 3$} and a generalization of the
  {K}auffman-{M}urasugi theorem}, Algebra i Analiz 2 (1990) 171--191

\bibitem{MR93c:57004}
\textbf{Yu\,V Drobotukhina}, \emph{Classification of projective {M}ontesinos
  links}, Algebra i Analiz 3 (1991) 118--130

\bibitem{MR1296890}
\textbf{Yu\,V Drobotukhina}, \emph{Classification of links in {${\bf R}{\rm
  P}\sp 3$} with at most six crossings [{MR} 93b:57006]}, from: ``Topology of
  manifolds and varieties'', Adv. Soviet Math. 18, Amer. Math. Soc.,
  Providence, RI (1994)  87--121

\bibitem{MR88h:05034}
\textbf{Jonathan~L Gross}, \textbf{Thomas~W Tucker}, \emph{Topological graph
  theory}, Wiley-Interscience Series in Discrete Mathematics and Optimization,
  John Wiley \& Sons Inc., New York (1987), a Wiley-Interscience Publication

\bibitem{MR29:621}
\textbf{J\,F\,P Hudson}, \textbf{E\,C Zeeman}, \emph{On combinatorial isotopy},
  Inst. Hautes \'Etudes Sci. Publ. Math.  (1964) 69--94

\bibitem{MR86e:57006}
\textbf{Vaughan F\,R Jones}, \emph{A polynomial invariant for knots via von
  {N}eumann algebras}, Bull. Amer. Math. Soc. (N.S.) 12 (1985) 103--111

\bibitem{MR1914297}
\textbf{Naoko Kamada}, \emph{On the {J}ones polynomials of checkerboard
  colorable virtual links}, Osaka J. Math. 39 (2002) 325--333

\bibitem{MR2001h:57007}
\textbf{Naoko Kamada}, \textbf{Seiichi Kamada}, \emph{Abstract link diagrams
  and virtual knots}, J. Knot Theory Ramifications 9 (2000) 93--106

\bibitem{MR88f:57006}
\textbf{Louis~H Kauffman}, \emph{State models and the {J}ones polynomial},
  Topology 26 (1987) 395--407

\bibitem{MR2000i:57011}
\textbf{Louis~H Kauffman}, \emph{Virtual knot theory}, European J. Combin. 20
  (1999) 663--690

\bibitem{MR1997331}
\textbf{Greg Kuperberg}, \emph{What is a virtual link?}, Algebr. Geom. Topol. 3
  (2003) 587--591 (electronic), math.GT/0208039

\bibitem{MR1999636}
\textbf{Maciej Mroczkowski}, \emph{Diagrammatic unknotting of knots and links
  in the projective space}, J. Knot Theory Ramifications 12 (2003) 637--651

\bibitem{MR2128052}
\textbf{Maciej Mroczkowski}, \emph{Polynomial invariants of links in the
  projective space}, Fund. Math. 184 (2004) 223--267

\bibitem{MR2002c:57009}
\textbf{Sam Nelson}, \emph{Unknotting virtual knots with {G}auss diagram
  forbidden moves}, J. Knot Theory Ramifications 10 (2001) 931--935

\bibitem{MR2002m:57011}
\textbf{Daniel~S Silver}, \textbf{Susan~G Williams}, \emph{Virtual knot
  groups}, from: ``Knots in Hellas '98 (Delphi)'', Ser. Knots Everything 24,
  World Sci. Publishing, River Edge, NJ (2000)  440--451

\end{thebibliography}
\def\cprime{$'$}

\Addresses

\end{document}